\documentclass[12pt]{article}

\usepackage{amsmath,amssymb,amsthm,amscd}
\usepackage[enableskew]{youngtab}
\usepackage[usenames]{color} 
\usepackage{young}
\usepackage{pstricks,pst-node,pst-tree}
\usepackage{ytableau}
\usepackage{hyperref}

\begin{document}

\newcommand{\oneb}{\overline{1}}
\newcommand{\twob}{\overline{2}}
\newcommand{\threeb}{\overline{3}}
\newcommand{\fourb}{\overline{4}}
\newcommand{\fiveb}{\overline{5}}
\newcommand{\jb}{\overline{j}}

\newcommand{\zerop}{0'}
\newcommand{\onep}{1'}
\newcommand{\twop}{2'}
\newcommand{\threep}{3'}

\newcommand{\threem}{-3'}
\newcommand{\twom}{-2'}
\newcommand{\onem}{-1'}
\newcommand{\fourm}{-4'}

\newcommand{\boldone}{\mathbf{1}}
\newcommand{\boldtwo}{\mathbf{2}}
\newcommand{\boldthree}{\mathbf{3}}
\newcommand{\boldfour}{\mathbf{4}}
\newcommand{\boldthreep}{\mathbf{3'}}

\newcommand{\boldoneb}{\mathbf{\overline{1}}}
\newcommand{\boldtwob}{\mathbf{\overline{2}}}
\newcommand{\boldthreeb}{\mathbf{\overline{3}}}
\newcommand{\boldfourb}{\mathbf{\overline{4}}}

\newcommand{\four}{4}
\newcommand{\three}{3}

\newcommand{\ybd}{\Yboxdim25pt}

\newcommand{\End}{{\rm{End}\ts}}
\newcommand{\Hom}{{\rm{Hom}}}
\newcommand{\ch}{{\rm{ch}\ts}}
\newcommand{\non}{\nonumber}
\newcommand{\wh}{\widehat}
\newcommand{\ot}{\otimes}
\newcommand{\la}{\lambda}
\newcommand{\La}{\Lambda}
\newcommand{\al}{\alpha}
\newcommand{\be}{\beta}
\newcommand{\ga}{\gamma}
\newcommand{\si}{\sigma}
\newcommand{\vp}{\varphi}
\newcommand{\de}{\delta^{}}
\newcommand{\om}{\omega^{}}
\newcommand{\hra}{\hookrightarrow}
\newcommand{\ve}{\varepsilon}
\newcommand{\ts}{\,}
\newcommand{\qin}{q^{-1}}
\newcommand{\tss}{\hspace{1pt}}
\newcommand{\U}{ {\rm U}}
\newcommand{\Y}{ {\rm Y}}
\newcommand{\CC}{\mathbb{C}\tss}
\newcommand{\SSb}{\mathbb{S}\tss}
\newcommand{\ZZ}{\mathbb{Z}\tss}
\newcommand{\NN}{\mathbb{N}\tss}
\newcommand{\AAA}{\mathbb{A}\tss}
\newcommand{\Z}{{\rm Z}}
\newcommand{\Ac}{\mathcal{A}}
\newcommand{\Pc}{\mathcal{P}}
\newcommand{\Qc}{\mathcal{Q}}
\newcommand{\Sc}{\mathcal{S}}
\newcommand{\Bc}{\mathcal{B}}
\newcommand{\Ec}{\mathcal{E}}
\newcommand{\Hc}{\mathcal{H}}
\newcommand{\Ar}{{\rm A}}
\newcommand{\Ir}{{\rm I}}
\newcommand{\Zr}{{\rm Z}}
\newcommand{\gl}{\mathfrak{gl}}
\newcommand{\Pf}{{\rm Pf}}
\newcommand{\oa}{\mathfrak{o}}
\newcommand{\spa}{\mathfrak{sp}}
\newcommand{\g}{\mathfrak{g}}
\newcommand{\ka}{\kappa}
\newcommand{\p}{\mathfrak{p}}
\newcommand{\sll}{\mathfrak{sl}}
\newcommand{\agot}{\mathfrak{a}}
\newcommand{\qdet}{ {\rm qdet}\ts}
\newcommand{\sdet}{ {\rm sdet}\ts}
\newcommand{\Gr}{ {\rm Gr}\tss}
\newcommand{\sgn}{ {\rm sgn}\ts}
\newcommand{\wt}{{\rm wt}\ts}
\newcommand{\Sym}{\mathfrak S}
\newcommand{\fand}{\quad\text{and}\quad}
\newcommand{\Fand}{\qquad\text{and}\qquad}
\newcommand{\vt}{{\tss|\hspace{-1.5pt}|\tss}}
\renewcommand{\theequation}{\arabic{equation}}

\newtheorem{thm}{Theorem}[section]
\newtheorem{lem}[thm]{Lemma}
\newtheorem{prop}[thm]{Proposition}
\newtheorem{cor}[thm]{Corollary}
\newtheorem{conj}[thm]{Conjecture}

\theoremstyle{definition}
\newtheorem{defin}[thm]{Definition}

\theoremstyle{remark}
\newtheorem{remark}[thm]{Remark}
\newtheorem{example}[thm]{Example}

\newcommand{\bth}{\begin{thm}}
\renewcommand{\eth}{\end{thm}}
\newcommand{\bpr}{\begin{prop}}
\newcommand{\epr}{\end{prop}}
\newcommand{\ble}{\begin{lem}}
\newcommand{\ele}{\end{lem}}
\newcommand{\bco}{\begin{cor}}
\newcommand{\eco}{\end{cor}}
\newcommand{\bde}{\begin{defin}}
\newcommand{\ede}{\end{defin}}
\newcommand{\bex}{\begin{example}}
\newcommand{\eex}{\end{example}}
\newcommand{\bre}{\begin{remark}}
\newcommand{\ere}{\end{remark}}
\newcommand{\bcj}{\begin{conj}}
\newcommand{\ecj}{\end{conj}}

\newcommand{\bal}{\begin{aligned}}
\newcommand{\eal}{\end{aligned}}
\newcommand{\beq}{\begin{equation}}
\newcommand{\eeq}{\end{equation}}
\newcommand{\ben}{\begin{equation$}}
\newcommand{\een}{\end{equation$}}

\newcommand{\bpf}{\begin{proof}}
\newcommand{\epf}{\end{proof}}

\newcommand{\aaa}{\mathbf{a}\tss}
\newcommand{\klmn}{K^\nu_{\nu\mu}}
\newcommand{\clmn}{c^\nu_{\la\mu}}
\newcommand{\cpmn}{c^\nu_{(p)\mu}}
\newcommand{\kkmn}{K^\nu_{\kappa\mu}}
\newcommand{\ckmn}{c^\nu_{\kappa\mu}}
\newcommand{\Trunc}{\textrm{Tr}}
\newcommand{\eval}{\textnormal{ev}}
\newcommand{\row}{\textrm{row}}
\newcommand{\col}{\textrm{col}}
\newcommand{\orow}{\overline{\textrm{row}}}
\newcommand{\ocol}{\overline{\textrm{col}}}
\def\ol#1{\overline{#1}}
\newcommand{\rhoz}{\rho^{(0)}}
\newcommand{\rhoo}{\rho^{(1)}}
\newcommand{\rhot}{\rho^{(2)}}
\newcommand{\rhoth}{\rho^{(3)}}
\newcommand{\rhol}{\rho^{(l)}}
\newcommand{\rhon}{\rho^{(n)}}
\newcommand{\rhom}{\rho^{(m)}}
\newcommand{\Tbar}{\overline{T}} 
\newcommand{\Sbar}{\overline{S}}
\newcommand{\one}{\mathbf{1}}
\newcommand{\zero}{\mathbf{0}}
\def\oneover#1#2{\frac{1}{c_{(1)#1}^{#1}-c_{(1)#2}^{#2}}}
\def\twoover#1#2{c_{(1)#1}^{#1}-c_{(1)#2}^{#2}}
\def\arho#1{\rho^{(#1)}}
\def\rbar#1{\hspace{2pt}|\overline{r_{#1}}|\hspace{2pt}}
\def\bar#1{\hspace{2pt}|\overline{#1}|\hspace{2pt}}
\def\K#1#2#3{K_{#1 #2}^{#3}}
\def\C#1#2#3{C_{#1 #2}^{#3}}
\def\Khat#1#2#3{\hat{K}_{#1 #2}^{#3}}
\def\asym#1#2{#1^{(#2)}}
\def\au#1#2{a_{#1(#2)}}
\newcommand{\cb}{\color{blue}}
\newcommand{\yvc}{\Yvcentermath1}
\renewcommand{\leq}{\leqslant}
\renewcommand{\geq}{\geqslant}
\renewcommand{\tilde}{\widetilde}
\renewcommand{\succeq}{\succcurlyeq}
\renewcommand{\preceq}{\preccurlyeq}
\def\barredtableaux#1#2#3#4#5#6
{\setlength{\unitlength}{0.75em}
\begin{center}
\begin{picture}(19,19.2)

\put(2,0){\line(0,1){19}}
\put(12,0){\line(0,1){4}}
\put(2,0){\line(1,0){10}}
\put(2,19){\line(1,0){17}}
\put(12,4){\line(1,0){2}}
\put(14,4){\line(0,1){4}}
\put(17,8){\line(0,1){6}}
\put(19,14){\line(0,1){5}}

\put(2,4){\line(1,0){8}}
\put(8,2){\line(1,0){4}}
\put(10,4){\line(0,-1){2}}
\put(8,2){\line(0,1){2}}
\put(8.5,2.7){#2}
\put(4.5,1.7){#1}
\put(2,10){\line(1,0){3}}
\put(3,8){\line(1,0){14}}
\put(5,10){\line(0,-1){2}}
\put(3,8){\line(0,1){2}}
\put(3.5,8.7){#4}
\put(6.5,5.7){#3}
\put(9.5,11){$\vdots$}
\put(2,16){\line(1,0){7}}
\put(7,14){\line(1,0){12}}
\put(9,16){\line(0,-1){2}}
\put(7,14){\line(0,1){2}}
\put(7.5,14.7){#5}
\put(6.5,17.2){#6}
\end{picture}
\end{center}
\setlength{\unitlength}{1pt}
}


\def\beql#1{\begin{equation}\label{#1}}

\title{\Large\bf Raising operators and the Littlewood--Richardson polynomials}

\author{{\sc Alex Fun}\\[10mm]
School of Mathematics and Statistics\\
University of Sydney,
NSW 2006, Australia\\
{\tt alex.fun\hspace{0.09em}@\hspace{0.1em}sydney.edu.au}
}

\date{} 


\maketitle

\vspace{15 mm}



%

\begin{abstract}
\noindent We use Young's raising operators to derive a Pieri rule for the ring generated by the indeterminates $h_{r,s}$ given in Macdonald's 9th Variation of the Schur functions. Under an appropriate specialisation of $h_{r,s}$, we derive the Pieri rule for the ring $\La(a)$ of double symmetric functions, which has a basis consisting of the double Schur functions. Together with a suitable interpretation of the Jacobi--Trudi identity, our Pieri rule allows us to obtain a new proof of a rule to calculate the Littlewood--Richardson polynomials, which gives a multiplication rule for the double Schur functions.
\end{abstract}

\setcounter{tocdepth}{2}
\tableofcontents

\section{Introduction}
The ring of symmetric functions $\La$ has a distinguished basis consisting of the Schur functions $s_\la$, parametrised by all partitions $\la$. The {Littlewood--Richardson rule} calculates the coefficients $\clmn$ occuring in the expansion of the product of two Schur functions:
\[
s_\la s_\mu = \sum_\nu \clmn s_\nu,
\]
summed over all partitions $\nu$. The coefficients $\clmn$ are nonnegative integers which play an important role in combinatorics \cite{Macdonald}, representation theory \cite{Sagan}, and geometry \cite{Fulton:yt}. There are now many versions of the Littlewood--Richardson rules, each utilising different combinatorial objects, which all give the same result; the survey paper by van Leeuwin \cite{Vanleeuwen} describes the equivalences between different versions of the Littlewood--Richardson rules.

Let $a = (a_i), i \in \ZZ$ be a sequence of variables. A generalisation of the ring $\La$ is the ring double symmetric functions $\La(a)$,  which is the ring of symmetric functions in $x_1, x_2, \dots$ with coefficients in $\ZZ[a]$. The ring $\La(a)$ has a basis consisting of the double Schur functions $s_\la(x|\!|a)$, parametrised by all partitions, see \cite{Molev:lrp} and \cite[Remark 2.11]{Okounkov:newton}. These are a multiparameter generalisation of the classical Schur functions. The classical ring of symmetric functions $\La$ is recovered by specialising $a$ to the sequence of zeroes. The Littlewood--Richardson polynomials $\clmn(a)$ arise as the structure coefficients in the following expansion:
\[
s_\la(x|\!|a)s_\mu(x|\!|a) = \sum_\nu \clmn(a) s_\nu(x|\!|a),
\]
summed over partitions $\nu$. A summary of the applications of the polynomials $\clmn(a)$ in combinatorics, geometry and resentation theory can be found in \cite{Molev:lrp}: under certain specialisations of the sequence $a$, the polynomials $\clmn(a)$ arise from a multiplication rule for equivariant Schubert classes \cite{Tao}, see also recent work  in the theory of isotropic Grassmanians \cite{bkt1, bkt2} and affine Grassmanians \cite{Shimozono1, Shimozono2}. Moreover, under a different specialisation of $a_i$ the polynomials $\clmn(a)$ give a multiplication rule for virtual quantum immanants and higher Capelli operators \cite{Okounkov}, \cite{Okounkov:qi}. Furthermore, after a shift of variables, the Littlewood--Richardsonn polynomials become the structure coefficients for the symmetric polynomials in the basis of the generalised factorial Schur functions, first calculated in \cite{Molevsagan}.

Compared to the Littlewood--Richardson coefficients, there are fewer rules to calculate the polynomials $\clmn(a)$. The first Graham positive \cite{Graham} rule  was given in \cite{Tao}, from the context of equivariant Schubert calculus, whereas an earlier rule given in \cite{Molevsagan} lacks the positive property. The rule in \cite{Tao} was expressed using the combinatorics of puzzles, see \cite{Justin} for another rule expressed using puzzles.

The first positive rule given in terms of tableaux were independently derived in \cite{Kreiman} and \cite{Molev:lrp}. Although the rules given are equivalent, the methods used to derive the rules in 
\cite{Kreiman} and \cite{Molev:lrp} are quite different. In \cite{Molev:lrp}, a recurrence relation of \cite{Molevsagan} is used, whereas \cite{Kreiman} generalises a concise proof of the classical Littlewood--Richardson rule by Stembridge \cite{Stembridge}, which relies on the definition of the Schur polynomial as a ratio of alternants.

The aim of this paper is to provide another method to calculate $\clmn(a)$ in terms of tableaux. The tableaux used here to calculate $\clmn(a)$ are the same as the ones in \cite{Molev:lrp}; instead, the novelty of our approach lies in the method by which $\clmn(a)$ will be calculated. See Remark \ref{rmk:alternativerule} for a possible way of using skew tableaux of shape $\nu/\mu$ to express $\clmn(a)$ using our approach. There are two main results. The first is the introduction of a Pieri rule for the ring $A$ associated with the indeterminates $h_{r,s}$ introduced in the 9th Variation of Macdonald's \cite{Macdonald:tv}. We do this by generalising a method of Tamvakis \cite{Tamvakis} involving Young's raising operators \cite{Young}, and show that these raising operators may be used on the polynomials $h_{r,s}$. Under an appropriate specialisation of the $h_{r,s}$ we derive a Pieri rule 
for the ring of double symmetric functions, which lets us compute $h_p (x|\!|a)s_\la(x|\!|a)$, the product  of a double complete symmetric function $h_p (x|\!|a)$ and a general double Schur function $s_\la(x|\!|a)$. The second main result is a method of calculating $\clmn(a)$ using our Pieri rule and the Jacobi--Trudi identity. The Jacobi--Trudi identity allows us to expand the Schur function in terms of an alternating sum of products of double complete symmetric functions. We use our Pieri rule, and then cancel out unwanted summands from the alternating sum. This is a generalisation of the method of using a Pieri rule and the Jacobi--Trudi identity to calculate the coefficients   $\clmn$ in the classical case; see Gasharov \cite{Gasharov}, Remmel and Shimozono \cite{RS}, and also Tamvakis \cite{Tamvakis}. 

Under a different specialisation, the 9th Variation polynomials $h_{r,s}$ \cite{Macdonald:tv} leads to the generalised Frobenius-Schur functions of \cite{ORV}. An area of interest is to calculate the structure coefficients arising from these functions, which would help explain the comultiplication structure in the ring $\La(a)$ of double symmetric functions; see \cite{Molev:dsci}. This in turn can be related to the equivariant cohomology of infinite grassmanians; see recent work by Liou and Schwarz \cite{LiouSchwarz}.

\section{Raising operators}

In this section we introduce Young's raising operators \cite{Young} and how they are applicable to the polynomials $h_{r,s}$ defined as the 9th variation in Macdonald \cite{Macdonald:tv}. We start with some preliminary definitions which are commonly found in the literature.

\bde
Let $\alpha$ be an \emph{integer sequence}, which is the sequence $(\alpha_1, \alpha_2, \dots)$ of integers. The sum of a pair of integers sequences $\alpha$ and $\beta$ is $\alpha + \beta = (\alpha_1 + \beta_1, \alpha_2+ \beta_2, \dots)$. An integer sequence $\alpha$ is called an \emph{integer vector} if it only has a finite number of nonzero entries. Suppose there exist a postive integer $l$, such that $\alpha_l$ is the right most nonzero entry. Then, $l$ is called the \emph{length} of $\alpha$, which is denoted  by $l(\alpha)$. If $\al$ is the sequence of zeroes, then the length of $\alpha$ is set to be 0. We will identify an integer vector $\alpha$ with the finite sequence, $(\alpha_1, \alpha_2, \dots, \alpha_l)$. For an integer sequence $\alpha$ and $k >0$, we define a \emph{truncation} of $\alpha$ to be the integer vector $\alpha^k = (\alpha_1, \dots,\alpha_k)$. We say that an integer vector $\alpha$ \emph{contains} another integer vector $\beta$ if $\alpha_i \geq \beta_i$ for each $i$, and denote this by $\alpha \supseteq \beta$.  
  A \emph{composition} is an integer vector such that every entry is nonnegative; a \emph{partition} is a composition whose entries weakly decrease, reading left to right. 

Let $\al$ be a composition and $l = l(\al)$. The \emph{diagram} of $\al$ is the finite collection of boxes, left justified, with rows numbered 1 to $l$, starting from the top and ending at the bottom, such that row $i$ has $\al_i$ boxes, $i = 1,\dots, l$. We will identify $\al$ with its diagram. A box $\ga$ in $\al$ is denoted $\ga \in \al$, and if this box is in row $i$, column $j$ of $\al$, then we write $\ga = (i,j)$. The \emph{content} of $\ga$ is $c(\ga) = j-i$. 

The \emph{skew diagram} $\la/\mu$ corresponding to a pair of partitions $(\la,\mu)$ with $\la \supseteq \mu$ is the set of boxes in $\la$ but not in $\mu$. A \emph{horizontal strip} is a skew diagram such that there is at most one box in each column.

\ede

In the 9th Variation of Macdonald's \cite{Macdonald:tv}, a commutative ring $A$ is defined as the ring generated by the independent indeterminates $h_{r,s}$, for $(r \geq 1, s \in \ZZ)$, over $\ZZ$. For convenience, define $h_{0,s} = 1$ and $h_{r,s} = 0$ for all $r <0$ and all $s \in \ZZ$. Define an automorphism $\tau$ of the ring $A$ generated by the $h_{r,s}$ by $\tau(h_{r,s}) = h_{r,s+1}$ for all $r,s$.  Let $\mu$ be an integer vector with $l = l(\mu)$ and $\beta$ an integer sequence. We will define elements $h_{\mu,\beta}$ of $A$, corresponding to the pair $(\mu,\beta)$, to be
\beq\label{eqn:x99}
h_{\mu,\beta} = h_{\mu_1,\beta_1}h_{\mu_2,\beta_2}\dots h_{\mu_l,\beta_l},
\eeq
written in this order.


We define \emph{raising operators} $R_{st}$, $1\leq s <t $, acting on integer sequences $\sigma$ by raising the $s$-th component of $\sigma$ by 1 and decreasing the $t$-th component of $\sigma$ by 1. That is, $R_{st} \sigma = (\dots, \sigma_s +1,\dots, \sigma_t -1,\dots)$. These operators were first introduced by Young \cite{Young}. An alternate explanation of how Young used these operators is given by Garsia \cite{Garsia}. Let $R$ be a monomial in the raising operators $R_{st}$. We will also call $R$ a raising operator and say that the operator $R$ acts on the integer vector, integer sequence pair $(\mu,\beta)$ via $R(\mu,\beta) = (R\mu,R\beta)$. We wish to write the following alternating sum:
\beq
s_{\mu,\beta} = \prod_{1\leq s < t \leq l} (1 - R_{st}) h_{\mu,\beta},\label{eqn:schur}
\eeq
which we interpret to mean the following. Let $R$ be a monomial in the $R_{st}$'s occuring in the expansion of the product in equation (\ref{eqn:schur}). Then we let the polynomial $R h_{\mu,\beta}$ be the polynomial $h_{R\mu,R\beta}$, respecting the order in which the factors of $h_{\mu,\beta} = h_{\mu_1,\beta_1}\dots h_{\mu_l,\beta_l}$ are written in equation (\ref{eqn:x99}).  In other words, the action of the operator $R$ on $h_{\mu,\beta}$ is induced by the action of $R$ on the pair $({\mu,\beta})$.  We will call $s_{\mu,\beta}$ the \emph{Schur polynomial} corresponding to the pair $(\mu,\beta)$.

We explain why we require equation (\ref{eqn:schur}). For convenience, let $R_l = \displaystyle\prod_{1 \leq i < j \leq l}(1 -R_{ij})$. We show that the Schur polynomial can be given by the \emph{Jacobi--Trudi identity}:
\begin{eqnarray}
s_{\mu,\beta} &=& R_l \left(h_{\mu,\beta}\right)\nonumber\\
&=&\det(h_{\mu_i+j-i,\beta_i +j-i})_{1\leq i,j\leq l}.\label{eqn:jtrudi}
\end{eqnarray}

Thus, the way we interpret equation (\ref{eqn:schur}) is that it is a way of symbolically expanding  the determinant in expression (\ref{eqn:jtrudi}), such that there might be more terms than necessary; for example, let $\mu$ be a composition with $l(\mu)>1$, then  $R_{12}R_{23} h_{\mu,\beta} = R_{13}h_{\mu,\beta}$, but these occur with opposite signs in the expansion of $s_{\mu,\be}$ in equation (\ref{eqn:schur}).

When $\beta = (0,1,2, \dots)$ and $\mu$ is a partition, equation (\ref{eqn:jtrudi})  is essentially equation (9.1') of Macdonald \cite{Macdonald:tv}. We now show that equation (\ref{eqn:jtrudi}) is true by suitably modifying an argument of Tamvakis \cite{Tamvakis}. Consider the ring $\mathbb{B}$ of Laurent polynomials in the (noncommuting) variables $x_{i,k}$, for $i \in \ZZ, k = 1,2,\dots$, with coefficients in $\ZZ$. For an integer vector $\mu$ with $l = l(\mu)$, and integer sequence $\beta$, we let
\[x^\mu_\beta = x^{\mu_1}_{\beta_1,1}x^{\mu_2}_{\beta_2,2}\dots x^{\mu_l}_{\beta_l,l}\] be a monomial. Then the raising operator $R$ acts on this monomial by $R x^\mu_\beta = x^{R\mu}_\beta$. Note that the second subscript $j$ in each $x_{i,j}$ makes the action of $R$ on $x^\mu_\beta$ ordered in the same way as with the action of $R$ on $h_{\mu,\beta}$. Let $\psi_\mu: \mathbb{B} \to \AAA$ be the $\ZZ$-linear map which takes $x^m_{nk}$ to $h_{m,m+n-\mu_k}$. Thus, we consider $Rh_{\mu,\beta}$ as the  image of $x^{R\mu}_\beta$ under the map $\psi_\mu$.

Let $R = R_{ij}$ for a pair $(i,j)$ with $1\leq i<j \leq l$. Consider the action of $R$ on $x^\mu_\beta$. It is equivalent to multiplying $x^\mu_\beta$ by  $x_{\beta_i,i}x^{-1}_{\beta_j,j}$, thus we have:
\beq\label{eqn:det}
\begin{split}
R_{l(\mu)} x^\mu_\beta &= \prod_{1\leq i< j\leq l}({1-R_{ij}}) x^\mu_\beta\\
		&= \prod_{1\leq i< j\leq l}{(1-x_{\beta_i,i}x^{-1}_{\beta_j,j})} x^\mu_\beta\\
		&= \det (x^{\mu_i-i+j}_{\beta_i,i})_{1\leqslant i,j \leqslant l},
\end{split}\eeq
with the last line following from the Vandermonde identity:
\[
\prod_{1\leqslant i < j \leqslant l} (x_{\beta_j,j} - x_{\beta_i,i}) = \det(x_{\beta_i,i}^{j-1})_{1\leqslant i,j\leqslant l}.
\]
Now apply $\psi_\mu$ to both sides of equation (\ref{eqn:det}) and equation (\ref{eqn:jtrudi})  is proven.

For each pair of integers $(i,j)$, $1\leq i < j$ define the operator $\overline{R}_{ij}$ acting on an integer sequence $\alpha$ as follows:
\[\overline{R}_{ij}\alpha = (\alpha_1, \dots, \alpha_{i-1}, \alpha_j-1, \alpha_{i}+1, \dots, \alpha_{j-1}, \alpha_{i+1}, \alpha_{j+1},\dots).\]
$\overline{R}_{ij}$ swaps the entries in the $i$-th and $j$-th place of $\alpha$ and then decreases the $i$-th entry by 1, and increases the $j$-th entry by 1. Note that $\ol{R}_{ij}$ is not the inverse of $R_{ij}$; in fact, it is equivalent to applying $R_{ij}$ and then swapping entries in the $i$-th and $j$-th place. Furthermore, $\ol{R}_{ij}$ is an involution on the set of integer sequences.

\bpr[Straightening law]
Suppose we have $(\mu,\beta)$, which is a pair of a integer vector and integer sequence respectively. Let  $(\mu',\beta')= (\overline{R}_{i,i+1}\mu, \overline{R}_{i,i+1}\beta)$, for a $i < l(\mu)$. Then, $s_{\mu,\beta} = - s_{\mu',\beta'}$. 
\epr

\bpf
The Schur polynomial $s_{\mu,\beta} = \det A$ for some matrix $A$ as defined from equation (\ref{eqn:jtrudi}). Similarly, $s_{\mu',\beta'} = \det B$ for some matrix $B$. Then we may obtain $B$ from $A$ by swapping row $i$ and $i+1$ of $A$. 
\epf


For the rest of this section, fix an integer sequence $\beta$. Let $\la$ and $\mu$ be compositions. From now on we define $h_\la$ to be the polynomial  $h_\la = h_{\la,\beta}$, and $s_\mu$ to be the polynomial $s_\mu = s_{\mu,\beta}$. Suppose $\la = (0, \dots, 0, p)$, a composition with length $e$ and let $l = l(\mu)$.  We wish to compute $h_\la s_{\mu}$, the product equal to $h_{p,\beta_e}s_{\mu,\beta}$ and write this product as a sum of the Schur polynomials. Recall the truncated integer vector $\beta^f = (\beta_1, \beta_2, \dots, \beta_f)$, for all $f > 0 $. Denote by $(\mu,p)$ and  $(\beta^l,\beta_{e})$ the concatenations $(\mu,p) = (\mu_1,\dots,\mu_l,p)$ and $(\beta^l,\beta_{e}) = (\beta_1,\dots,\beta_l,\beta_{e})$. We have that
\beq\begin{split}
h_{p,\beta_e} s_\mu &= h_{p,\beta_e} R_{l} h_\mu\\
&= R_{l} \{h_\mu h_{p,\beta_e}\}\\
&=R_{l+1}\prod_{1\leq i \leq l}{(1+R_{i,l+1}+R^2_{i,l+1}+R_{i,l+1}^3+\dots)} \{h_{(\mu,p),(\beta^l,\beta_e)} \}\\
&=R_{l+1}\sum_R{R} h_{(\mu,p),(\beta^l,\beta_e)},\label{eqn:raising}
\end{split}\eeq
summed over monomials $R$ in the $R_{i,l+1}$'s, for all $1 \leq i \leq l$; the fourth equality holds because  \[(1+R_{i,l+1}+R^2_{i,l+1}+\dots)(1-R_{i,l+1}) = 1,\] since the action of $R_{i,l+1}$ on $h_{ij}$, $1 \leq i \leq l$, is nilpotent for all $j$. This leads to the Pieri rule.

\bpr[Pieri rule]
\label{prop:pieri}
Let $\mu$ be a partition with $l = l(\mu)$, let $p$ be a positive integer and $e$ be an integer. Let $\beta' = (\beta^l,\beta_e-p)$. Then 
\[h_{p,\beta_e}s_\mu = \sum_{\sigma} s_{\mu+\sigma,\beta'+\sigma},\]
summed over compositions $\sigma$ such that $\sigma$ has $p$ boxes and has length at most $l+1$.
\epr

\newcommand{\ds}{\displaystyle}
\bpf
From equation \ref{eqn:raising} we have that $h_{p,\beta_e} s_\mu = R_{l+1}\ds\sum_R{R} h_{(\mu,p),(\beta^l,\beta_e)}$, summed over monomials $R$ in the $R_{i,l+1}$'s, for all $1 \leq i \leq l$. For such a $R$, we have that
\[
R h_{(\mu,p),(\beta^l,\beta_e)} =  h_{\mu+\sigma,\beta'+\sigma},
\]
for some $\sigma$ such that $\sigma$ has $p$ boxes and has length at most $l+1$. That is, $R$
 acts on both the compositions  $(\mu,p)$ and $(\beta^l,\beta_e)$ by removing up to $p$ boxes from row $l+1$ of  each composition and adding them to the previous rows. The proposition follows.
\epf

\section{The double symmetric functions}

The ring of double symmetric functions is defined by Molev \cite{Molev:lrp} as follows. Let $a = (a_i)$, $i\in \ZZ$ be a sequence of variables. Consider the ring of polynomials
$\ZZ[a]$ in the variables $a_i$ with integer coefficients. Consider another infinite set of
variables $x = (x_1, x_2, \dots )$ and for each nonnegative integer $n$ denote by $\Lambda_n(a)$ the ring
of symmetric polynomials in $x_1,\dots , x_n$ with coefficients in $\ZZ[a]$. When $a$ is specialised to the sequence of zeroes, this reduces to the usual ring of symmetric polynomials, see e.g. Macdonald \cite{Macdonald}. Many of the properties we will discuss for $\Lambda_n(a)$ are analogues of similar properties of the ring of symmetric polynomials. The ring $\Lambda_n(a)$ is filtered by the usual degrees of polynomials in $x_1, \dots, x_n$ with the $a_i$ considered to have the zero degree. The evaluation map
\beq\label{eqn:inverse}
\varphi_n : \Lambda_n(a) \to \Lambda_{n-1}(a), \text{\hspace{1cm} $P(x_1,\dots , x_n) \to P(x_1,\dots , x_{n-1}, a_n)$},
\eeq
is a homomorphism of filtered rings so that we can define the inverse limit ring $\Lambda(a)$ by
\[\Lambda(a)  = \varprojlim \Lambda_n(a), \text{\hspace{1cm} $n \to \infty$}\]
where the limit is taken with respect to the homomorphisms $\varphi_n$ in the category of
filtered rings. The ring $\Lambda(a)$ is equipped with the automorphism $\tau$ which takes $a_i$ to $a_{i+1}$, for all $i\in \ZZ$. 
\newcommand{\xa}{(x|\!|a)}
$\Lambda_n(a)$ is generated over $\ZZ[a]$ by \emph{double complete symmetric polynomials} $h_p\xa$, for all $p \geq 1$:
\[
h_p(x_1,\dots, x_n|\!|a) = \sum_{n\geq i_1 \geq i_2 \geq \dots \geq i_p\geq 1} (x_{i_1} - a_{i_1})(x_{i_2} - a_{i_2-1})\dots(x_{i_p} - a_{i_p-p+1}),
\]
which have a stability property \cite{Okounkov}; i.e. they are compatible with respect to the homomorphisms (\ref{eqn:inverse}):
\beq\label{eqn:stability}
\varphi_n: h_p(x_1,\dots, x_n|\!|a) \mapsto h_p(x_1,\dots, x_{n-1} |\!|a).
\eeq
Thus, the ring $\Lambda(a)$ is generated by the {\em double complete symmetric functions} $h_p\xa$:
\[
h_p(x|\!|a) = \sum_{i_1 \geq \dots \geq i_p\geq 1} (x_{i_1} - a_{i_1})(x_{i_2} - a_{i_1-1})\dots(x_{i_p} - a_{i_p-p+1}).
\]

There are distinguished bases of $\Lambda_n(a)$ over $\ZZ[a]$, 
one of which are the double Schur polynomials $s_\la(x_1,\dots,x_n|\!|a)$, parameterized
by all partitions $\la$.  
%
%
Let $\one = (0,1,2,\dots)$. We will define the double Schur functions by fixing $\beta = \one$ and choosing a specialisation of the $h_{r,s}$ from the previous section.  For any $n >0$, we will set $h^{(n)}_{r,s} = \tau^s h_r(x_1,\dots,x_n|\!| a)$. This means the ring $A$ from the previous section specialises to $\Lambda_n(a)$, for a choice of $n >0$, when we set $h_{r,s} = h^{(n)}_{r,s}$. We have that the automorphism $\tau$ is the same as the operation on the sequence $a$ such that $\tau a_i = a_{i+1}$ for all $i \in \ZZ$.

For all $r \geq 1$, we have the following relation between $h_{r,s}$ and $h_{r,s-1}$:
\beq
h^{(n)}_{r,s} = h^{(n)}_{r,s-1} + (a_{s-r+1} - a_{n+s}) h^{(n)}_{r-1,s-1}.\label{eqn:relation}
\eeq
which can be directly calculated, or inferred from Molev \cite[Lemma 2.4]{Molev:lrp}. 

Let $\la$ be a partition with $l = l(\la)$. Then we define the \emph{double Schur polynomials} $s_\la(x_1,\dots,x_n|\!|a)$ by specialising the Jacobi--Trudi identity (equation \ref{eqn:jtrudi}):  
\begin{eqnarray}
s_\la(x_1,\dots,x_n|\!|a) &=& \det\left(h_{\la_i+i-j,j-1}^{(n)}\right)_{1\leq i,j \leq l}\label{eqn:jt}\\
&=& R_l h^{(n)}_{\la,\one^l}\nonumber
\end{eqnarray}
where here we mean that $h^{(n)}_{\la,\one^l} = h^{(n)}_{\la_1,0} \dots h^{(n)}_{\la_l,l-1}$. Recall that the last equality holds because we may use raising operators to rewrite the determinant. Since the $h^{(n)}_{r,s}$ are stable under (\ref{eqn:stability}), the polynomial  $s_\la(x_1,\dots,x_n|\!|a)$ is stable under the homormorphisms (\ref{eqn:inverse}). Thus, we may define the \emph{double Schur functions} $s_{\la}(x|\!|a) = (s_\la(x_1|\!|a),s_\la(x_1, x_2|\!|a), \dots) \in \Lambda(a)$, parametrized by all partitions $\la$, to be the basis of the inverse limit ring $\Lambda(a)$. 

For partitions $\la$ and $\mu$, the \emph{Littlewood--Richardson polynomials} $\clmn(a;n)$ are defined as the structure coefficients in the following expansion
\beq\label{eqn:stablelr}
s_\la(x_1,\dots,x_n|\!|a) s_\mu(x_1,\dots,x_n|\!|a) = \sum_{\nu} \clmn(a;n) s_\nu(x_1,\dots,x_n|\!|a),
\eeq
summed over partitions $\nu$. Furthermore, we define $\clmn(a)$ to be the structure coefficients in the expansion
\[
s_\la(x|\!|a) s_\mu(x|\!|a) = \sum_{\nu} \clmn(a) s_\nu(x|\!|a),
\]
summed over partitions $\nu$, and also call these Littlewood--Richardson polynomials, a terminology first introduced in \cite{Molev:lrp}.

The polynomials $\clmn(a)$ and $\clmn(a;n)$ are polynomials in $\ZZ[a]$, and in the case of $\clmn(a;n)$, are dependent on $n$.  However since the polynomial $s_\la(x_1,\dots,x_n|\!|a)$ is stable under the homormorphisms (\ref{eqn:inverse}), we have that $\clmn(a;n)$ does not depend on $n$ when $n$ is big enough. This is a remarkable fact which means that for such a sufficiently big $n$ the coefficient $\clmn(a) = \clmn(a;n)$. Our main aim is to use the Jacobi--Trudi identity (\ref{eqn:jt}) to calculate $\clmn(a)$. We first sketch a way in which $\clmn(a;n)$ may be calculated. First expand the determinant in equation (\ref{eqn:jt}) to obtain alternating summands consisting of products of double complete symmetric polynomials. Then calculate the product of each of these summands with  $s_\mu(x_1,\dots,x_n|\!|a)$, and  using the Pieri rule, decompose the result in terms of the double Schur polynomials. Hence each alternating summand from the expansion of (\ref{eqn:jt}) contributes to the polynomial $\clmn(a;n)$.

In the classical case such a calculation for the Littlewood--Richardson coefficients has been explored; see \cite{Gasharov, RS}. These authors used a `sign reversing involution' to simplify the contribution of the alternating summands appearing in equation (\ref{eqn:jt}) to the Littlewood--Richardson coefficient. We will use this idea to calculate $\clmn(a)$, with the details to follow in Section 3. This method motivates the following defintion of $\tau$, an automorphism of $\Lambda(a)$, which applies to $h_p(x|\!|a) \in \Lambda(a)$ in the following way:
\beq \label{eqn:BIG}
\tau h_p(x|\!|a) = h_p(x|\!|a) + a_{2-p} h_{p-1}(x|\!|a),
\eeq
for all $p \geq 1$.  With this definition of $\tau$, the Jacobi--Trudi identity (\ref{eqn:jt}) may be reinterpreted for $\Lambda(a)$:
\beq\begin{split}\label{eqn:JT}
s_\la(x|\!|a) &= \det\left(\tau^{j-1}h_{\la_i+i-j}(x|\!|a)\right)_{1\leq i,j \leq l}\\
&= R_l h_{\la,\one}.
\end{split}\eeq
Relations (\ref{eqn:relation})  and (\ref{eqn:BIG}) are essentially the same: The difference between these two relations is that we choose to ignore the coefficient depending on $n$ in (\ref{eqn:relation}). We do this because for the purposes of using equation (\ref{eqn:JT}) to calculate $\clmn$, since for a big enough $n$, the coefficient $\clmn(a) = \clmn(a;n)$ and no longer depends on $n$, so all contributions of a coefficient involving $n$ in $\clmn(a;n)$ must vanish as $n$ becomes very big. 

\subsection{Pieri rule: after specialisation}

From now onwards, we will work with the double Schur functions and double complete symmetric functions, and use relation (\ref{eqn:BIG}) in our  calculations; for the rest of this paper let $h_{r,s}$ denote $\tau^s h_r(x|\!|a)$ and $s_\la = s_\la\xa$, given by equation (\ref{eqn:JT}). 

Let $\mu$ be a partition and let $l = l(\mu)$. Recall that $h_{\mu,\be} = h_{\mu_1,\be_1}\dots h_{\mu_l,\be_l}$.  Let $\one'$ denote the concatenation  $(\one^l,e-p) = (0, 1,\dots, l-1,e-p)$. From Proposition \ref{prop:pieri}, we have the corresponding Pieri rule for $\Lambda(a)$:
\[\begin{split}h_{p,e}s_\mu &= \sum_{\sigma} s_{\mu+\sigma,\one'+\sigma}\\
&=\sum_{\sigma} R_{l+1} h_{\mu+\sigma,\one'+\sigma}
\end{split}\]
summed over compositions $\sigma$ with $p$ boxes, and of at most length $l+1$. We use  relation (\ref{eqn:BIG}) to rewrite
\beq\label{eqn:pieribig}
\sum_\sigma h_{\mu+\sigma,\one'+\sigma} = \sum_{\sigma'} d_{\sigma'}(a) h_{\mu+\sigma',\one^{l+1}},
\eeq
summed over compositions $\sigma' \subseteq \sigma$, for all $\sigma$ with $p$ boxes and of at most length $l+1$. Here, the coefficient $d_{\sigma'}(a)$ is a polynomial in $\ZZ[a]$. Thus,
\begin{eqnarray}
h_{p,e}s_\mu &=&  R_{l+1} \displaystyle\sum_{\sigma'} d_{\sigma'}(a) h_{\mu+\sigma',\one^{l+1}}  \nonumber\\
&=&\displaystyle\sum_{\sigma'} d_{\sigma'}(a) s_{\mu+\sigma',\one^{l+1}}\label{eqn:dsigma}
\end{eqnarray}

We claim that
\beq
\sum_{\sigma'} d_{\sigma'}(a) s_{\mu+\sigma',\one^{l+1}} 
= \sum_{\nu} \clmn(a) s_\nu(x|\!|a) \label{eqn:pierirule}
\eeq
summed over partitions $\nu = \mu+\sigma'$, such that $\nu/\mu$ is a horizontal strip, and $\la = (0,\dots,0,p)$, a composition of length $e+1$ (recall then that $h_{\la} = h_{p,e}$). Call $\sigma'$ a \emph{good} composition if $\nu/\mu$ is a horizontal strip, and \emph{bad} otherwise. Then, we claim that $d_{\sigma'}(a) = \clmn(a)$ when  $\sigma'$ is good; furthermore, if $\sigma'$ is bad then $\clmn(a) = 0$ and the contributions of $d_{\sigma'}(a) s_{\mu+\sigma',\one^{l+1}}$ for all bad $\sigma'$'s will cancel on the left hand side of equation (\ref{eqn:pierirule}). The aim of the rest of this section is to prove these claims, and provide a way of calculating $d_{\sigma'}(a)$. 

We first calculate $d_{\sigma'}(a)$. It will be useful to think of our compositions as diagrams. Letting $l = l(\mu)$, for a $1\leq j \leq l$, observe that if $\sigma_j$ contains $m$ boxes then $m$ boxes are added to row $j$ of $\mu$ when we form $\mu+\sigma$. Then, we will have the pair $(\mu_j + m ,	 j-1 + m)$ appearing as the $j$-th entry in the pair $(\mu+\sigma, \one^l + \sigma)$. This $j$-th entry corresponds to the polynomial $h_{\mu_j + m , j-1 + m}$. To obtain the right hand side of equation (\ref{eqn:pieribig}) we wish to reduce the second index down to $j-1$. Using relation (\ref{eqn:BIG}),
\begin{eqnarray}
h_{\mu_j + m , j-1 + m} &=& h_{\mu_j + m , j-1 + m -1} + a_{j-\mu_j} h_{\mu_j + m-1 , j-1 + m -1}\nonumber\\
&=& h_{\mu_j + m , j-1 } + \sum_{d=0}^{m-1} a_{j-m -\mu_j +d+1} h_{\mu_j + m-1 , j-1 + d}\label{eqn:relns},
\end{eqnarray}
and doing the same to $h_{\mu_j + m-1 , j-1 + d}$, and so on, we have \[
h_{\mu_j + m , j-1 + m} = \displaystyle\sum_{d = 0}^{m} K_d(a) h_{\mu_j + m - d , j-1},\]
where $K_d(a)$ is the following degree $d$ polynomial in $\ZZ[a]$:
\beq\label{eqn:kd}
K_d(a) = \sum_{b_1,\dots,b_d}\prod_{i=1}^d a_{b_i - \mu_j - m +  1}, 
\eeq
summed over integers $b_i$ such that $ j +d \leq b_d \leq \dots \leq b_1 \leq j + m-1	$.

For the case where $\sigma_{l+1}$ contains $p-k$ boxes, for a $0 \leq k \leq p$, we have the pair $(p-k, e-k)$ as the $l+1$-th entry in the pair $(\mu+\sigma, \one^l + \sigma)$. Since $e < l$, we have that $e - k < l$. We want this second index to be $l$. We may use relation (\ref{eqn:BIG}) to write
\[\begin{split}
h_{p-k,e-k} 
= h_{p-k,l} - \sum_{d = e-k}^{l-1} a_{d-p+k+2}h_{p-k-1, d},
\end{split}\]
and doing the same to $h_{p-k-1,d}$ and so on, we can write 
\beq\label{eqn:gd}
h_{p-k,e-k} = \sum_{d=0}^{p-k} G_d(a) h_{p-k-d,l}
\eeq
where $G_d(a)$ is the following  polynomial in $\ZZ[a]$ with degree $d$
\[G_d(a) =  \sum_{b_1,\dots,b_d}\prod_{i=1}^d (-1)^d a_{b_i-p+k+i+1},\]
summed over integers $b_i$ such that $e-k\leq b_1\leq \dots \leq b_d \leq l-1$.

We will now eliminate bad compositions from the sum on the left hand side of equation (\ref{eqn:pierirule}), using tableaux which we define in the next section. 

\subsection{Using tableaux to calculate $d_{\sigma'}(a)$}

We deduce from our calculations for $K_d(a)$ and $G_d(a)$ that if $\sigma \supseteq \sigma'$, then $h_{\nu,\one+\sigma}$ will contribute to a summand in $d_{\sigma'}(a)$. We will make the calculation of $d_{\sigma'}$ more precise using tableaux, which we now define. An example of all the following definitions exists at the end of this section.

Fix a composition $\la$. A \emph{reverse $\la$-tableau} $T$ of shape $\la$ is obtained by filling each box of $\la$ with a positive integer $k$ which is unbarred, or a positive integer $\overline{k}$ which is barred; further, in each row of $T$, the entries weakly decrease from left to right. Note that we do not impose any conditions on the columns of $T$. The following definitions are all \emph{associated} to the tableau $T$. If $\alpha = (i,j)$ is a box of $\la$, we call $T(\alpha) = T(i,j)$ the \emph{entry} of $T$ in box $\alpha$. The \emph{content} of box $\alpha$ is $c(\alpha) = c(i,j) = j-i$. 

We introduce the \emph{row order} on the boxes of $\la$ by reading the boxes in rows from bottom to top, from left to right of each row. We extend this ordering to the entries of a tableau $T$ of shape $\la$. Let $\alpha$ and $\beta$ be two boxes of $\la$. Then $T(\alpha)$ is \emph{before} $T(\beta)$ with respect to the row order if the box $\alpha$ is before $\beta$ with respect to the row order.

Similarly, introduce the \emph{column order} on the boxes of $\la$ by reading the boxes in columns left to right, from the bottom to the top of each column. Similarly, we will say $T(\alpha)$ is \emph{before} $T(\beta)$ with respect to the column order if the box $\alpha$ is before $\beta$ with respect to the column order.

Let $\mu$ be a composition and let $S = s_1s_2\dots s_t$, be a sequence of positive integers. We \emph{apply} $S$ to $\mu$ by forming a sequence of compositions from $\mu$ which terminates in $\nu$ in the following way:
\[
\mu = \rho^{(0)}\overset{s_1}{\to} \rho^{(1)} \overset{s_2}{\to} \dots \overset{s_r}{\to}  \rho^{(t)} = \nu,
\]
such that $\rho^{(i)}$, $i= 1, \dots, t$, are compositions and $\rho^{(i-1)}\overset{s_i}{\to} \rho^{(i)}$ means adding a box to the end of row $s_i$ of $\rho^{(i-1)}$ to form $\rho^{(i)}$. We say that $S$ \emph{takes} $\mu$ to $\nu$ (or $\nu$ is \emph{created} from $\mu$ using $S$), and denote this as $S:\mu \to \nu$.  We will say that $S$ is \emph{Yamanouchi} when applied to $\mu$ if $\rho^{(i)}$ is a partition for all $0 \leq i \leq t$, and \emph{not Yamanouchi} otherwise.

If $T$ is a reverse $\la$-tableau, we define the \emph{row} and \emph{column} word of $T$. The \emph{row word} $S^r$ \emph{corresponding} to $T$  is the sequence of barred entries in $T$ listed left to right, from the first barred entry to the last with respect to the row order. Similarly, the \emph{column word} $S^c$ \emph {corresponding} to $T$ is the sequence of barred entries in $T$ listed left to right, from the first to the last with respect to the column order. When writing the row and column words corresponding to $T$ we will omit the bars. For every $\alpha = (i,j) \in \la$, let $S^r(\alpha) = S^r(i,j)$ be the subsequence of $S^r$ consisting of the barred entries in $T$ listed up to, and including box $\alpha$, with respect to the row order. Let $S^c(\alpha) =S^c(i,j)$ be the subsequence of $S$ consisting of the barred entries in $T$ listed up to, and including box $\alpha$, with respect to the column order. We will let $\rho^r(\alpha) = \rho^r(i,j)$ (resp. $\rho^c(\alpha) = \rho^c(i,j)$) be the composition created from $\mu$ using $S^r(\alpha)$ (resp. $S^c(\alpha)$).

For each  $\alpha \in \la$, let $\rho(\alpha)$ be a composition. Then we define the \emph{weight} of an entry $T(\alpha)$ to be
\[
\eval(T(\alpha)) = a_{T(\alpha) - \rho(\alpha)_{T(\alpha)}} - a_{T(\alpha) -c(\alpha)}.
\]
The \emph{weight} of the tableau $T$ is the weight of all unbarred entries of $T$ multiplied together, denoted:
\[
\eval{(T)} = \prod_{\substack{\alpha \in \la \\ \alpha \text{ unbarred}}} (a_{T(\alpha) - \rho(\alpha)_{T(\alpha)}} - a_{T(\alpha) -c(\alpha)}).
\]

\bex
Let $\la = (3,5,2)$, which corresponds to the diagram
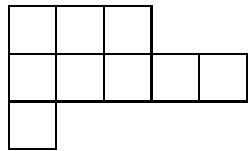
\[
\ydiagram{3,5,1}
\]
In turn, we have the row and column ordering on the boxes of $\la$, which we illustrate in the following two diagrams by filling in the boxes of $\la$ with integers so that the first box with respect to the ordering is labelled `1' and so on
\[\begin{array}{ccc}
\ytableaushort{789,23456,1} &\hspace{1cm}& \ytableaushort{357,24689,1}\\
\text{row ordering} & \hspace{1cm}& \text{column ordering}
\end{array}\]

Consider the following reverse $\la$-tableau $T$ of shape $\la$:
\[T = \ytableausetup{centertableaux}
\ytableaushort{
2{\ol{2}}{\ol{1}},
22{\ol{2}}11{\ol{1}},
1
}
\]
Then, the row word $S^r = 2121$ of $T$ is the sequence of barred integers of $T$, listed with respect to the row order, and the column word $S^c = 2211$ of $T$ is the sequence of barred integers of $T$, listed with respect to the column order. Let $\mu = (2,1)$ be a partition. Then the row word $S^r$ of $T$ takes $\mu$ to $\nu = (4,3)$ via the following sequence of compositions:
\[
\mu = \ydiagram{2,1} \overset{2}{\to} \ydiagram{2,2} \overset{1}{\to} \ydiagram{3,2} \overset{2}{\to} \ydiagram{3,3} \overset{1}{\to} \ydiagram{4,3} = \nu
\]
Since all of these are partitions, we say that $S^r$ is Yamanouchi when applied to $\mu$. For the box $(1,2)$ we have the subword $S^r(1,1) = 212$, the sequence of barred integers of $T$ listed up to box $(2,1)$ with respect to the row order. Correspondingly, the partition $\rho^r(1,2)$ is the partition $(3,3)$. We leave it to the reader to check that $S^c$ takes $(2,1)$ to $(4,3)$ but is not Yamanouchi.

For the purposes of the weight of the entry $T(2,4)$, let $\rho(\al) = \rho^c(\al)$. Then, the weight of $T(2,4)$ is
\[\begin{split}
\eval(T(2,4)) &= a_{1-\rho^c(2,4)_1} - a_{1-c(2,4)}\\
&=a_{-2} - a_{-1}.
\end{split}\]
\qed\eex

\subsection{Weights of tableaux express $d_{\sigma'}(a)$}
 
For the rest of this section let $\la = (0,\dots,0,p)$ of length $e+1$ and assume that $e$ is strictly less than the length of the partition $\mu$.  We will use reverse $\la$-tableaux to make sense of equation (\ref{eqn:pieribig}) and the coeffiecients $d_{\sigma'}(a)$ appearing in it. Recall that $d_{\sigma'}(a)$ is the coefficient of $h_{\mu+\sigma',\one^{l+1}}$ when we rewrite (\ref{eqn:pieribig}). We claim that
\beq\label{eqn:claim}
d_{\sigma'}(a) = \sum_{T} \eval(T),
\eeq
summed over all reverse $\la$-tableaux $T$ such that $T$ has row word $S^r: \mu \to \mu + \sigma'$. For example, suppose $\sigma' = \sigma$ is a composition with $p$ boxes, and is at most length $l+1$. Then, there is a unique reverse $\la$-tableau $T$ with row word $S^r$ such that $S^r:\mu \to \mu+\sigma$, so our claim in this instance is that $d_\sigma = 1$. We will prove the claim (\ref{eqn:claim}) for all $\sigma' \subseteq \sigma$ in the following paragraphs.

Let $T$ be a reverse $\la$-tableau, with row word $S^r:\mu \to \mu+\sigma'$, for a $\sigma' \subseteq \sigma$. Let us expand the weight of $T$ in terms of monomials. We create a tableau $U$, \emph{derived} from $T$, by doing the following: for every unbarred entry $x \in T$ we will either leave it unbarred, or put a prime on it, that is, we replace the entry $x$ with $x'$. The tableau $U$ inherits the definitions associated with $T$: for example, $S^r$, $\rho^r(\al)$ etc. The weight of an unbarred entry  $U(\alpha)$ is then taken to be $a_{T(\alpha) - \rho^r(\alpha)_{T(\alpha)}}$, and the weight of a unbarred primed entry $U(\alpha)$ is taken to be $-a_{|T(\alpha)|-c(\alpha)}$, where $|T(\alpha)|$ means disregard the prime on the entry in box $\alpha$. Thus we see that the weight of $T$ may be expressed as
\[\eval(T) = \sum_U \prod_{\substack{ \alpha \in \la,\\ \alpha \text{ unbarred, unprimed}}}a_{T(\alpha) - \rho^r(\alpha)_{T(\alpha)}}\prod_{\substack{\alpha \in \la,\\ \alpha  \text{ unbarred, primed}}}-a_{|T(\alpha)|-c(\alpha)},\]
summed over all tableaux $U$ derived from $T$.

\bre
Although we do not call them such, these tableaux $U$ introduced here are equivalent to the barred reverse $\la$-supertableaux introduced in Section 4 of Molev \cite{Molev:lrp}.
\ere


We now prove that each monomial appearing in the coefficient $d_{\sigma'}(a)$ (\ref{eqn:dsigma}) may be given by the weight of a tableau $U$ derived from a reverse $\la$-tableau $T$ with row word $S^r:\mu \to \mu+\sigma' = \nu$. 

We first show that we need only consider \emph{good} tableaux $U$, which are tableau $U$ that satisfy the following two conditions:
\begin{description}
\item[Cond. 1:] The maximum unbarred, unprimed entry appearing in $U$ is $l$.
\item[Cond. 2:] Let $y$ be the number of barred $l+1$'s appearing in $U$. Then there does not exist a primed entry $U(\alpha)$, for some $\alpha \in \la$, such that $|U(\alpha)| - c(\al) > l-y$. 
\end{description}
Any tableau $U$ which is not good is \emph{bad}. Note that since $S^r:\mu \to \nu$ and $l(\nu) \leq l+1$ we have that the maximum barred integer appearing is $l+1$. Let $\mathbb{U}$ denote the set of bad tableaux $U$. We construct a weight reversing involution on $\mathbb{U}$, such that a bad tableau $U$ in $\mathbb{U}$ is paired to a tableau $\tilde{U}$ in $\mathbb{U}$ with reverse weight to $U$. If $U$ is a bad tableau, we will call an entry $U(\be)$ \emph{bad} if $U(\be) > l$ is unprimed, or $U(\be)$ is primed and $|U(\be)| - c(\be) > l - y$; i.e. the entry $U(\be)$ violates Cond. 1 or Cond. 2 charecterising good tableaux. Then, let $\alpha \in \la$ be the box satisfying both of the following conditions:
\begin{description}
\item[Cond. 1:] $U(\al)$ is bad and the subscript of the weight of $U(\alpha)$ is the maximal subscript appearing in the weight of any bad entry of $U$. Call this subscript $k$.
\item[Cond. 2:] If there is more than one entry with weight equal to $a_k$ or $-a_k$: let $\alpha$ be the box containing the primed entry with weight equal to $-a_k$ if it exists, otherwise, $\alpha$ is the leftmost box containing a unprimed, unbarred entry with weight equal to $a_k$.  
\end{description}
Let $j$ be the column number of the box $\al$, so that $\al = (e+1,j)$. There are two cases to consider in the construction of $\tilde{U}$:


\emph{Case 1:} Suppose that $U(e+1,j)$ is unbarred and unprimed. Then, there exist a unique pair: a box $(e+1,j')$, with $j' \geq j$ and a positive integer $m$ such that $ m \geq |U(e+1,j'+1)|$ and $m -c(e+1,j') = k$. First, we argue that such an $m$ exists. Let $m(j') = k + c(e+1,j')$ for all $j' \geq j$. Then, this is a strictly increasing sequence of integers, since $c(e+1,j')$ strictly increases as $j'$ increases. On the other hand $|U(e+1,j')|$ weakly decreases as $j'$ increases. Second, $m$ is positive since $m - c(e+1,j) \geq k$ and  $e+1 \leq l(\mu)$. Third, the pair is unique since the subscripts  of the weights of any primed entries in $U$ strictly decrease, reading left to right along the row.  To form $\tilde{U}$, we will remove the entry $U(e+1,j)$, and move all entries from box $(e+1,j+1)$ to $(e+1,j')$ inclusive, one box to the left. Now insert the entry $m'$ into box $(e+1,j')$. Note that there are no primed entries between box $(e+1,j)$ and $(e+1,j')$, by assumption of maximality of $k$. Thus by construction, the tableau $\tilde{U}$ formed this way has opposite weight to $U$.

\emph{Case 2:} Suppose that $U(e+1,j)$ is primed and that $k = l+1 - y'$, for some $0 \leq y' \leq y$. Then there exist a box $(e+1,j')$, for some minimal $j' \leq j$, such that there are $y'$ barred $l+1$'s strictly to the left of it. We will remove the entry in box $(e+1,j)$, move all entries from box $(e+1,j')$ to $(e+1,j-1)$ inclusive one box to the right, and then insert an unbarred $l+1$ in box $(e+1,j')$. Suppose $k > l+1$; in this case we remove the entry in box $(e+1,j)$, move all entries from box $(e+1,1)$ to $(e+1,j-1)$ inclusive, one box to the right, and insert an unbarred $k$ in box $(e+1,1)$. Again by construction, the tableaux formed in either of these ways has opposite weight to $U$.

We now give an example which illustrates this involution.

\bex
Let $\la = (0,0,5)$, $\mu = (2,1)$, and $\nu = (2,1,3)$. Then consider the following tableaux $U$:
\[
U = \ytableaushort{ {\ol{3}} {\ol{3}} 3 {\ol{3}} 2}
\]
The entries $U(3,3)$, and $U(3,5)$ both have weight $a_1$. Since they are both unprimed, this is dealt with in Case 1. Thus, we pick $\al$ to be the leftmost box, i.e. $\al = (3,3)$. We have the tableau $\tilde{U}$: 
\[
\tilde{U} = \ytableaushort{ {\ol{3}} {\ol{3}}  {\ol{3}} {2'} 2}
\]
which is formed from $U$ by removing the entry $U(3,3)$, moving the entry $U(3,4)$ one spot to the left, and inserting the entry $2'$ into the box $(3,4)$. This new entry has weight equal to $-a_{2 - c(3,4)} = -a_1$.  The entry $2'$ and the box $(3,4)$ are a unique pair in the sense that there is no other box of $\la$ and a primed entry that can go in that box and have weight equal to $-a_1$.

Note that $\tilde{U}$ has a primed entry which has weight $-a_1$, and `1' is in fact the maximum subscript of the weight of any entry in $\tilde{U}$. Thus, this tableau is dealt with in Case 2. We restore $U$ by noting that box (3,3) is the leftmost box of $\tilde{U}$ which has two barred $3$'s strictly to the left of it.

\qed
\eex

We now show that the weight of monomial tableaux $U$ may be used to represent monomials appearing in $d_{\sigma'}$ (\ref{eqn:pieribig}). Let $U$ be a good tableau. Then for each $1 \leq j \leq l+1$,  we have that there are $\sigma'_j$ barred $j$'s in $U$. For any $1 \leq j\leq l$, we consider the following argument: let there be $d$ unbarred, unprimed $j$'s in $U$, and let the product of their weights be the monomial $V$, written left to right with respect to the row order; thus, the subscripts of $V$ weakly decrease. Let $\sigma_j = \sigma'_j + d$, and we will show that $V$ is equal to a monomial in $K_{d}$, the coefficient of $h_{\mu_j + \sigma'_j , j-1}$  in equation (\ref{eqn:pieribig}) formed by applying the relation (\ref{eqn:BIG}) to $h_{\mu_j + \sigma_j,j - 1+\sigma}$.  First note that $V$, like $K_{d}$, is a degree $d$ monomial in $\ZZ[a]$. From equation (\ref{eqn:kd}), the maximum subscript appearing in a factor of $V$ is $j - \mu_j$, while the minimum is $j- \mu_j - \sigma'_j$. Since the subscripts weakly decrease, reading left to right, the subscripts of $V$ agree with the subscripts of a unique monomial appearing in the expansion of $K_d(a)$ in equation (\ref{eqn:kd}). Thus, the weight of the subtableau of $U$ containing only $j$'s is equal to a unique monomial in $K_d(a)$.

 Abusing notation, let there be $d$ primed entries in $U$. Now consider the monomial $W$ equal to the product of  the weights of all primed entries in $U$, written left to right, with respect to the row order. Thus, the subscripts of $W$ strictly decrease, reading left to right. We will show that $W$ corresponds to a unique monomial appearing in the coefficient $G_d(a)$ in equation (\ref{eqn:gd}). We have that the maximum subscript appearing is $l-y$ where $y$ is the number of barred $l+1$'s, since $U$ is a good tableau. Noting that $y = p - k -d$ in equation (\ref{eqn:gd}), the integer $l-y$ agrees with the maximum subscript appearing in  $G_d(a)$. The minimum subscript appearing in $W$ is $1 - c(e+1,p) = 2 - p + e$. This agrees with the minimum subscript appearing in  $G_d(a)$. Thus, we see that the subscripts of $W$ agree with the subscripts of a unique monomial appearing in the expansion of $G_d(a)$.

Let $\nu = \mu + \sigma'$ and $\la = (0,\dots,0,p)$ of length $e+1 \leq l(\mu)$. Thus we conclude that $d_{\sigma'}(a)$ from equation (\ref{eqn:dsigma}) may be expressed using the weights of reverse $\la$-tableaux;
\[
d_{\sigma'} = \sum_{T} (a_{T(\alpha) - \rho^r(\alpha)_{T(\alpha)}} - a_{T(\alpha) -c(\alpha)}),
\]
summed over tableaux $T$ of shape $\la$, with row word $S^r: \mu \to \nu$, such that $\nu = \mu + \sigma'$. The main result of this section is the following proposition, which is a restatement of the claim made in (\ref{eqn:pierirule}).

\subsection{Pieri rule: Statement and Proof}
\bpr\label{prop:doublepieri}
Let $\mu$, $\nu$ be partitions, and $\la = (0,\dots,0,p)$, of length $e+1$, such that $1 \leq e < l-1$. If $\mu \not\subseteq \nu$, we have $\clmn(a) = 0$. Otherwise, 
\beq\label{eqn:doublepieri}
\clmn(a) = \sum_{T} (a_{T(\alpha) - \rho(\alpha)_{T(\alpha)}} - a_{T(\alpha) -c(\alpha)}),
\eeq
summed over reverse $\la$-tableaux $T$, with row word $S^r: \mu \to \nu$ that is Yamanouchi, and $\rho(\alpha)$ is defined by $S^r(\alpha): \mu \to \rho(\alpha)$.
\epr


\bpf
If $\mu \not\subseteq \nu$, we have $\clmn(a) = 0$, since there is no row word that can create $\nu$ from $\mu$. Recall $\sigma'$ is bad if $\nu/\mu$ is not a horizontal strip. If $\sigma'$ is bad, we will call the diagram $\nu = \mu + \sigma'$ bad as well.  We will show that if $S^r$ is not Yamanouchi, then the coefficient $\clmn(a) = 0$.  We do this by cancelling all contributions of $d_{\sigma'}(a) s_{\mu+\sigma',\one^{l+1}}$ on the left hand side of equation (\ref{eqn:pierirule}), for all bad $\sigma'$'s.

Let $\nu$ be bad, and $S^r$ the row word, which is weakly decreasing, reading left to right, which takes $\mu$ to $\nu$. Let $\mathbb{T}$ be the set of reverse $\la$-tableaux $T$ with row word equal to $S^r$. Let $i$ be minimal such that $\nu_{i+1} - \mu_{i} > 0$. Such an $i$ exists since $\nu$ is bad. Let $\tilde{\nu} = \ol{R}_{i,i+1} \nu$, and $\tilde{S}$ be the sequence of integers, weakly decreasing reading left to right, which takes $\mu$ to $\tilde{\nu}$. Then, the composition $\tilde{\nu}$ is bad. Let $\tilde{\sigma'}$ be the composition equal to difference $\tilde{\nu} - \mu$. Let $\tilde{\mathbb{T}}$ be the set of reverse $\la$-tableaux $T$ with row word equal to $\tilde{S}$.

We claim the following:
\[
\sum_{T \in \mathbb{T}} \eval(T) = \sum_{T \in \mathbb{\tilde{T}}} \eval(T).
\]
We will prove this by constructing a weight preserving bijection between monomial tableaux $U$ derivable from $T \in \mathbb{{T}}$ and monomial tableaux $\tilde{U}$ derivable from $\tilde{T} \in \mathbb{\tilde{T}}$. An example follows at the end.

Let $U$ be a monomial tableau derived from $T \in \mathbb{T}$. Let $d = \nu_{i+1}- \nu_i  - 1$.
We construct $\tilde{U}$ with equal weight to $U$; there are two cases depending on whether $d\geq 0$ or $d<0$:

\emph{Case 1:} Suppose $d \geq 0$. Let $X$ be the subtableau of $U$ containing $d$ barred $i+1$'s, counting left from the right most $i+1$ in $U$. We will replace the  $d$ barred $i+1$'s in $X$ with barred $i$'s to create a tableau with a new row word that will create $\tilde{\nu}$ from $\mu$. We do this in the following way. First, we do not do anything with the primed entries in $U$.  Let $X^{i+1}$ denote the sequence of unprimed $i+1$'s in $X$, reading left to right. Let $Y$ denote the sequence of nonprimed entries equal to $i$ in $U$, reading left to right. Delete the entries in $X$ which are not primed, and also delete all unprimed entries equal to $i$, thus creating empty boxes. Let $X^i$ be the sequence created from $X^{i+1}$ by replacing all the unbarred $i+1$'s in $X^{i+1}$  with unbarred $i$'s, and all the barred $i+1$'s in $X^{i+1}$ with barred $i$'s. We will now fill in the newly created empty boxes as follows. Insert entries in the empty boxes of our tableau, from left to right, with the entries from $Y$, read left to right, and then the entries from $X^i$, read left to right. 

Due to the previous processes applied to $U$, the entries of our current tableau, read left to right,  may no longer weakly decrease. So the next step is to fix the order in which the barred entries equal to $i$ and $i+1$ occur. Call a primed $i$ (resp. $i+1$) \emph{badly ordered} if it is to the left of a $i+1$ (resp. right of a $i$). We describe a process which fixes badly ordered primed $i$'s. A very similar process will fix the badly ordered primed $i+1$'s; see the example below. Starting with the rightmost badly ordered primed $i$, do the following: Delete the badly ordered primed $i$. Move the $i+1$ to the right of it one box to the left. Now insert a primed $i+1$ in the blank box. Keep repeating this process on the next rightmost badly ordered primed $i$ until none are left.

We end up with a tableau with entries that weakly decrease, read left to right, and this is the tableau $\tilde{U}$ paired with $U$.

\emph{Case 2:} Suppose $d < 0$. Then we can undo the processes described in Case 1. Let $X$ be the subtableau of $T$ containing $|d|$ barred $i$'s, counting left from the rightmost $i$ in $T$.  Then we may reverse the process described in Case 1.

Note that the involution induces a natural pairing between the entries of $U$ and those of $\tilde{U}$.  

We check that the weight of $U$ is the same as the weight of $\tilde{U}$, for $U$ a Case 1 tableau. Let the box $\alpha$ contain an entry of $U$ that was changed from an $i$ to $i+1$, or vice versa, in the creation of $\tilde{U}$.   Suppose $\alpha$ contains an unbarred $i+1$ that was changed into an unbarred $i$. Then by the definition of $d$, in $U$ there are $\nu_{i} +1$ barred $i+1$'s strictly to the left of $\alpha$. The unbarred entry $i$ in box $\be$ of $\tilde{U}$ which is paired to $U(\alpha)$ by the involution has  $\nu_{i}$ barred $i$'s strictly to the left of it. Thus, the weight $\eval(U(\alpha)) = i + 1 - \rho^r(\alpha)_{i+1} = i + 1 - (\nu_{i} +1) $ is equal to $\eval(\tilde{U}(\be)) = i - \rho^r(\alpha)_i = i  - \nu_{i}$. On the other hand, suppose $\alpha$ contained a badly ordered primed $i$ that was changed into a primed  $i+1$. Then $\tilde{U}(\be) = (i+1)'$ is the entry paired with $U(\al)$, and $\be$ is one box to the right of $\al$. Since $c(\beta) = c(\al)+1$, we conclude the weights of $U(\al)$ and $\tilde{U}(\be)$ are equal. A similar argument follows for badly ordered primed $i+1$'s that were changed into primed  $i$'s.


\epf

\bex
Let $\la = (8)$, $\mu = (3,2)$ and $\nu = (3,5)$, a bad diagram. Then, the composition $\tilde{\nu} = \ol{R}_{1,2}(\nu) = (4^2)$ is bad as well, since $\tilde{\nu}/\mu$ is not a horizontal strip. Let $U$ be the following monomial tableau
\[
U =\ytableaushort{ {\ol{2}}  {\ol{2}} {\ol{2}} {2'} 2 {1'} 1 {1'} }
\]
which has row word $S^r = 222$. We have that $d = \nu_{i+1} - \nu_i - 1 = 1$. Then, the tableau $X$ is the subtableau of entries from boxes (1,3) up to (1,5), since $X$ contains $d = 1$ barred 2's, counting left from the rightmost 2. The sequence of unprimed entries in $X$, read left to right, is $X^2 = \ol{2} 2$. The sequence $X^1 = \ol{1} 1$ is formed from $X^2$ by replacing all 2's by 1's. The sequence of unprimed 1's in $T$ is $Y = 1$. We now delete unprimed entries to form:
\[
\ytableaushort{ {\ol{2}}  {\ol{2}} {} {2'} {} {1'} {} {1'} }
\]
Now, fill in the blank boxes with entries from $Y$, then $X^1$, read left to right:
\[
\ytableaushort{ {\ol{2}}  {\ol{2}} 1 {2'} {\ol{1}} {1'} 1 {1'} }
\]
The entries do not weakly decrease, read left to right, so we must fix the badly ordered primed entries, by swapping the $2'$ with the 1 on its left, and then replacing $2'$ with $1'$. This forms the tableau $\tilde{U}$:
\[
\tilde{U} = \ytableaushort{ {\ol{2}}  {\ol{2}} {1'}  1 {\ol{1}} {1'} 1 {1'} }
\]
We claim the weight of $U(1,5)$ is equal to the weight of $\tilde{U}(1,7)$. This is because $2 - \rho^r(1,5)_2 = 2 - 5$ and $1 - \tilde{\rho}^r(1,7)_1 = 1-4$. We claim the weight of $U(1,4)$ is equal to the weight of $\tilde{U}(1,3)$. This is true since $c(1,3) = c(1,4) -1$, so $1 - c(1,3) = 2 - c(1,4)$.

\qed
\eex

\section{Littlewood--Richardson polynomials}

Let $\mu$, $\la$ be partitions, and $l = l(\mu)$, and recall the integer vector $\one = (0,1,\dots)$ and the Jacobi--Trudi identity (\ref{eqn:JT}):
\beq\begin{split}\label{eqn:JT2}
s_{\mu,\one} =& \det(h_{\mu_i+j-i,j-1})_{1\leq i,j\leq l}.
\end{split}\eeq
Let $\pi = (l-1, l-2, \dots, 0)$ and $\Sym_l$ denote the symmetric group on $l$ elements. For each $\omega \in \Sym_l$, define the composition $\la^\omega = \omega(\la+\pi_{l})-\pi_{l}$ and let $\sgn(\la^\omega) = \sgn(\omega)$, the parity of the permutation $\omega$. We may write
\[
s_{\mu,\one} = \sum_{\kappa} \sgn{(\kappa)} h_{\kappa, \one}
,\]
summed over $\kappa = \la^\omega$, for all $\omega\in \Sym_l$. This is just an expansion of the determinant (\ref{eqn:JT2}) into an alternating sum. For each $\kappa$, define $K_{\kappa\mu}^{\nu}(a)$ as the coefficients appearing in the expansion 
\[
h_{\kappa,\one} s_{\mu,\one} = \sum_{\nu} K_{\kappa\mu}^{\nu}(a) s_{\nu,\one}.
\]
Then, we have that
\beq\label{eqn:weaklr}
\clmn(a) = \sum_{\kappa} \sgn(\kappa) K_{\kappa\mu}^{\nu}(a).
\eeq
summed over $\kappa = \la^\omega$, for all $\omega\in \Sym_l$.
	
In the classical case, when $\kappa$ is a partition and $a$ is the sequence of zeroes, the coefficient $K_{\kappa\mu}^{\nu}(a)$ are just the \emph{Kostka numbers}. 
 If $\kappa$ is a diagram with only one row then this is the Pieri rule (\ref{eqn:doublepieri})  and thus $\kkmn = \ckmn$. The aim of the rest of this paper is to eliminate unwanted coefficients $\kkmn$ from the above alternating sum. First we give a formula for $\kkmn(a)$ using the Pieri rule. 

\bpr\label{prop:kkmn}
Let $\kappa = \la^\omega$, for some $\omega \in \Sym_l$. Then
\beq\label{eqn:kostka}
\kkmn(a) = \sum_T \prod_{\underset{\alpha \text{ unbarred}}{\alpha \in \la}} \left(a_{T(\alpha)-\rho^r(\alpha)_{T(\alpha)}} - a_{T(\alpha)-c(\alpha)}\right),\eeq
summed over reverse $\kappa$-tableaux $T$, such that each $T$ has row word $S^r: \mu \to \nu$ which is Yamanouchi.
\epr


\bpf
The proof follows from repeated applications of the Pieri rule (\ref{eqn:doublepieri}). Let $l' = l(\la)$.

We have that 
\begin{eqnarray}
h_\kappa s_\mu &=& h_{\kappa_1} \dots h_{\kappa_{l'}} s_\mu \nonumber\\
&=& h_{\kappa_1} (\dots(h_{\kappa_{l'-1}}(h_{\kappa_{l'}}s_\mu))\dots) \label{eqn:rsdt}
\end{eqnarray}
where we evaluate each multiplicative pair using the Pieri rule, starting with $h_{\kappa_{l'}}s_\mu$. Each multiplication produces a tableau of shape $(0,\dots,0,\ka_i)$, of length $i$. We stack these tableau on top of each other to form a tableau $T$ of shape $\ka$; that is, the $i$-th row of $T$ is equal to the tableau formed from the $i$-th multiplicative pair in expression \ref{eqn:rsdt}. Furthermore, this tableau must  contain a row word $S^r: \mu \to \nu$.
\epf

\bth\label{thm:lr}
Let $\la$, $\mu$, and $\nu$ be partitions. If $\nu \not\subseteq \mu$, the coefficient $\clmn(a) = 0$. If $\mu \subseteq \nu$, we have that
\beq\label{eqn:lr}
\clmn(a) = \sum_{T} \prod_{\substack{\alpha \in \la\\ T(\alpha) \text{ unbarred}}}(a_{T(\alpha)-\rho^r(\alpha)_{T(\alpha)}} - a_{T(\alpha) - c(\alpha)}),
\eeq
where the sum is taken over reverse $\la$-tableaux $T$ obeying the following. First, the column word $S^c$ of $T$ is Yamanouchi and $S^c:\mu \to \nu$. Secondly, the entries in $T$ strictly decrease down each column.
\eth


The statement that $\clmn(a) = 0$ if $\mu \not\subseteq \nu$ follows from the Pieri rule (Proposition \ref{prop:doublepieri}). We will split the proof up into sections but first we introduce some terminology.  We will call a tableau $T$ \emph{good} if it appears in the sum (\ref{eqn:lr}). On the other hand, a \emph{bad} tableau $T$ is one appearing in the expansion of $\kkmn$ for all $\ka = \la^\omega$ in equation (\ref{eqn:kostka}) which is not good. The following is an equivalent description of bad tableaux: 

Let $\kappa$ be the shape of a bad tableau $T$, and $l = l(\kappa)$. For each row $i$, let $T^{\geq i}$ denote the subtableau of $T$ consisting of entries in row $i$ and below. If a box $\alpha$ is in row $i$ of $\kappa$, let $L(\alpha)$ be the column word corresponding to the the barred entries of $T^{\geq i}$ in the boxes before and including $\alpha$, with respect to the column order.  Then a tableau $T$ is bad if and only if $T$ has one of the following properties:
\begin{itemize}
	\item[(P1)] There exist a row $i$ and a box $\alpha$ in row $i$ of $\kappa$ such that the sequence $L(\alpha)$ is not Yamanouchi when applied to $\mu$.
	\item[(P2)] There exist a row $i$ of $\kappa$ such that the subtableau of $T$ formed from rows $i$ and $i+1$ of $T$ is not column strict.
	\item[(P3)] There is a row $i$ of $\ka$ such that $\kappa_i < \kappa_{i+1}$.
\end{itemize}

We split the proof up  Theorem \ref{thm:lr} into two sections.   In the first section we will describe an involution on the set of bad tableaux. Namely, we pair a bad tableau $T$ to another bad tableau $\tilde{T}$ of shape $\tilde{\ka} = R_{i,i+1}\ka$  appearing in the sum (\ref{eqn:weaklr}), for some $1 \leq i \leq l$. This is not a weight preserving involution, however it is close to one, as we will discover later. In the second section, we will use a sequence of lemmas to show that it is possible to cancel out the weights of bad tableaux from the sum (\ref{eqn:weaklr}). An example which ties these two sections together will follow at the end.


If $T$ is a bad tableau of shape $\ka$, then $\kappa_m \neq \kappa_{m+1} - 1$ for all rows $m$ of $T$. This is because $\la_i \geq \la_{i+1}$ which necessarily means that $\ka_i \geq \ka_{i+1}$ or $\ka_i \leq \ka_{i+1} -2$.  Furthermore, there exist a unique pair of integers $(i,j)$, subject to both of the following conditions  on $T$:
\begin{itemize}
\item[(C1)] The row number $i$ is maximal such that one of properties (P1), (P2) or (P3) hold for $T$.
\item[(C2)] The column number $j$ is minimal so one of the following (mutually exclusive) conditions hold for $T$:
\begin{itemize}
	\item[(C2a)] The property (P1) holds for $T$ and $\alpha = (i,j)$, and $T(i,j) > T(i+1,j)$.
	\item[(C2b)] The property (P2) holds for $T$ and $\alpha = (i,j)$; by this we mean $T(i,j) \leq T(i+1,j)$. 
	\item[(C2c)] The property (P3) holds for $i$, and there is no column $j$ such that (P1) and (P2) hold for box $\alpha = (i,j)$ in row $i$ of $T$. 
\end{itemize}
\end{itemize}
In the case where condition (C2c) holds, let $j = \kappa_{i+1} + 1$. For all cases, we will call the tableau $T$ \emph{bad} in row $i$, column $j$, and we have that $\ka_{i+1} \geq k_i + 2$.

%

We begin constructing our involution on the set of bad tableaux. Throughout this construction, if $\al$ and $\be$ are boxes in a row of $\ka$, when we write ``between boxes $\al$ and $\be$'' we mean that the boxes $\al$ and $\be$ are included in this range. If we want to exclude either of these boxes from the range we will specifically say so.

Let $T$ be a bad tableau with shape $\kappa$. We will construct $\tilde{T}$, a bad tableau paired to $T$, such that: 1) The shape $\tilde{\kappa}$ of $T$ is equal to $\ol{R}_{i,i+1} \kappa$. 2) The row word $\tilde{S}^r$ corresponding to $T$ takes $\mu$ to $\nu$. 

To construct $\tilde{T}$ we use a sequence of processes, which are summed up by the following diagram:
\[
T \overset{\psi_1}{\longrightarrow} T^1 \overset{\psi_2}{\longrightarrow} \tilde{T},
\]
where we start with the tableau $T$, and apply processes $\psi_1$ and $\psi_2$ in the order indicated by the arrows. This will create a sequence of tableaux involving the intermediate tableau $T^1$, and this sequence terminates at $\tilde{T}$, which is the bad tableaux paired with $T$. The first process, $\psi_1$, is called a ``tail swap'' and the second, $\psi_2$, is called ``reorder barred entries''. The detailed description of each of these processes appear in the following subsections with the same name.

\subsection{$\psi_1$: tail swap}

In this process, we take the tableau $T$ which is bad in row $i$, column $j$ and create a tableau $T^1$, of shape $\tilde{\ka}$, which preserves the `bad in row $i$, column $j$' condition. We first describe the process for tableaux which obey conditions (C2b), (C2c), and leave the argument for (C2a) till last, and then we will examine the properties of $T^1$.

\emph{Suppose $T$ obeys condition (C2b):} Let $T(i,j) = b$ and $T(i+1,j) = c$, and from the  definition of $\alpha = (i,j)$ in condition (C2b) we have $b \leq c$. We call the pair $b$ and $c$ a \emph{bad column pair}. We do the following; let $X$ be the subtableau of entries in boxes $(i,j)$ to $(i,\kappa_i)$ of $T$, and $Y$ the subtableau of entries in boxes $(i+1,j+1)$ to $(i,\kappa_{i+1})$ of $T$. Form the tableau $T^1$ by \emph{swapping} the subtableaux $X$ and $Y$. Swapping means two things:
\begin{itemize}
	\item[1.] The shape of $T^1$ is $\tilde{\ka} = \ol{R}_{i,i+1}\ka$.
\item[2.] The entries of $T^1$ in boxes $(i+1,j+1)$ up to $(i+1, \ka_i+1)$ are equal to the entries of $X$, read left to right, and the entries of $T^1$ in boxes $(i,j)$ up to $(i, \ka_{i+1}-1)$ are equal to the entries of $Y$, read left to right.
\end{itemize}
There is a natural \emph{pairing} of boxes affected by the tail swap; e.g. if $\alpha = (m,n)$ is a box in $X$, then the box $\alpha' = (m+1,n+1)$ is the box paired with $\alpha$ by the tail swap. In a similar vein, the entry $T(\alpha)$ is paired with the entry $T^1(\alpha')$. This pairing holds for the other tail swaps described for cases (C2a) and (C2c) as well.

\emph{Suppose $T$ obeys condition (C2c):}
Let $Y$ be the subtableau of entries of $T$ in boxes $(i+1,\kappa_i+2)$ to $(i+1,\kappa_{i+1})$. Then move $Y$ to the end of row $i$ of $T$ to form $T^1$. We can think of the process for $T$ obeying condition (C2c) as a special case of the process for tableaux obeying condition (C2b); that is we have that the subtableau $X$ is empty. 

The case where $T$ obeys condition (C2a) involves an additional step, which we now describe.

\emph{Suppose $T$ obeys condition (C2a):} Let $T(i,j) = b$, and $T(i+1,j) = c$. We will also call the pair $b$ and $c$ a \emph{bad column pair}. By definition of $\al = (i,j)$ in (C2a), we have that the sequence  $L(\alpha)$ is not Yamanouchi when applied to $\mu$. Therefore, we have that $b = c+1$ and the entry $T(i,j)$  is a barred $c+1$.
\bde\label{def:q}
Define the column number $q$ to be maximal so that the entry  $T(i,q)$ is a $b$, barred or unbarred. 
\ede
Suppose that there are $s$ barred $b$'s in boxes $(i,j)$ up to $(i,q)$. We now examine the structure of the entries of the row $i$ and $i+1$ of $T$.  For any $k\geq 1$, a \emph{block} of unbarred $k$'s is an uninterrupted sequence, reading left to right, of unbarred $k$'s in a row of $T$. Between boxes $(i,j)$ and $(i,q)$, we have $s$ disjoint blocks of unbarred $b$'s, with each block to the right of a barred $b$ (starting with the one in box $(i,j)$). Note that some of these blocks may be empty. Let $x_i$, $i = 1, \dots, s$, be the number of entries in each block respectively, reading the blocks left to right. 
\bde\label{def:r} Define the column number $r$ to be minimal such that there are $s$ barred $c$'s between box $(i+1,j+1)$ and $(i+1,r)$.
\ede
This column number exists since the row word of $T$ is Yamanouchi. We have that $T(i+1,r) = \overline{c}$ by definition of $r$. Between boxes $(i+1,j+1)$ and $(i+1,r)$ we have $s$ disjoint blocks of unbarred $c$'s, with each block to the left of a barred $c$. Let $y_i$, $ i = 1 ,\dots, s$, be the number of entries in each block respectively, reading left to right. There are two cases depending on whether $q \geq r$ or $q < r$:

\emph{Case 1:} If $q \geq r$, we will form an intermediate tableau $T^{\frac{1}{2}}$ by doing the following process, which we call \emph{fixing the column ordering}: Replace the entries in boxes $(i,j)$ up to $(i,r-1)$ such that the new entries consist of $s$ blocks, each containing $y_i$, $1 \leq i \leq s$, unbarred $b$'s, with each block to the right of a barred $b$ (starting with the one in box $(i,j)$). Let $\sigma$ be the subdiagram of $\ka$ containing the boxes from box  $(i+1,j+1)$ up to $(i+1,r)$, then $(i, r)$ up to $(i,q)$. Then, starting from box $(i+1,j+1)$ and ending in box $(i,q)$, replace the entries in the boxes of $\sigma$ with $s$ blocks, each containing $x_i$ unbarred $c$'s, $1 \leq i \leq s$, such that each block is before a barred $c$, with respect to the row order imposed on $\sigma$. The tableau formed this way is $T^{\frac{1}{2}}$. Note that in the tableau $T^{\frac{1}{2}}$ we have that the entries in boxes $(i,j)$ up to $(i,r-1)$ consist solely of $b$'s, and the entries in the boxes of $\sigma$ consist solely of $c$'s. In particular $|T^{\frac{1}{2}}(i+1,r)| = |T^{\frac{1}{2}}(i,r)| = c$.

Let $X$ be the subtableau of $T^{\frac{1}{2}}$ containing the entries from box $(i,r)$ up to $(i,\kappa_i)$, and $Y$ the subtableau of $T^{\frac{1}{2}}$ containing entries from box $(i+1, r+1)$ up to $(i+1,\kappa_{i+1})$. Now swap $X$ with $Y$ to obtain the tableau $T^1$. Note that that the shape of $T^1$ is $\tilde{\kappa} = R_{i,i+1}\kappa$ and the rows of $T^1$ weakly decrease. 

\emph{Case 2:}
We now deal with the case where $q < r$. Define $X$ (resp. $Y$) to be the subtableau of entries of $T$ from box $(i,q+1)$ up to $(i, \kappa_i)$ (resp. $(i+1,q+2)$ up to $(i, \kappa_{i+1})$). Swap the subtableaux $X$ and $Y$ and the tableau obtained is defined to be $T^{\frac{1}{2}}$. Again, the shape of $T^\frac{1}{2}$ is $\tilde{\kappa} = R_{i,i+1}\kappa$ and the rows of $T^{\frac{1}{2}}$ weakly decrease.

If $q = r -1$, set $T^1 = T^{\frac{1}{2}}$. If $q < r - 1$, in the tableau $T^{\frac{1}{2}}$ we have that the entries in boxes $(i,j)$ up to $(i,q)$ consist solely of $b$'s. Let $\sigma$ be the subdiagram of $\tilde{\kappa}$ consisting of the boxes $(i,j)$ up to $(i,q+1)$, then $(i,q+1)$ up to $(i,r-1)$. Then, the entries in $\sigma$ consist solely of $c$'s. Recall that the tableau $T^{\frac{1}{2}}$ formed in Case 1 has a very similar property. Now, we apply a process which we will also call fixing the column ordering: Replace the entries of $T^{\frac{1}{2}}$ in boxes $(i,j)$ up to $(i,r-1)$ such that the new entries consist of $s$ blocks, each containing $y_i$, $1 \leq i \leq s$, unbarred $b$'s, with each block to the right of a barred $b$ (starting with the one in box $(i,j)$). Then, replace the entries in boxes $(i+1,j+1)$ up to $(i+1,q+1)$ with $s$ blocks, each containing $x_i$ unbarred $c$'s, $ i = 1, \dots, s$, with each block to the left of a barred $c$. This creates the tableau $T^1$.

We have just described a process of creating the tableau $T^1$ from a tableau $T$ obeying (C2a), (C2b), or (C2c). The following properties of $T^1$ will be relevant when we want to show the processes described are an involution on the set of bad tableaux:
\begin{itemize}
	\item[1.] $T^1$ is bad in row $i$, column $j$.
	\item[2.] The rows of $T^1$ weakly decrease, left to right.
	\item[3.] $T^1$ is of shape $\tilde{\ka} = \ol{R}_{i,i+1} \ka$.
	\item[4.] If $T$ obeys condition (C2a), and falls under Case 1 (with $q \geq r$), then $\tilde{T}$ will obey condition (2a), and will fall under Case 2 with $q < r-1$. The converse holds as well: if $T$ obeys condition (C2a), and falls under Case 2 with $q < r-1$, then $\tilde{T}$ will obey condition (2a), and will fall under Case 1 (with $q \geq r$).
\end{itemize}
We will now describe the second process, to be applied to $T^1$.

\subsection{$\psi_2$: reorder barred entries}
We require this process because the row word formed from the barred entries of $T^1$ might not be Yamanouchi when applied to $\mu$. To fix this we will rearrange barred entries equal to at most $c$ in rows $i$ and $i+1$ of $T^1$.

For each $2 \leq k \leq c$, let $r_k$ be the number of barred $k$'s in every row below, and including row $i$. For each $1\leq k < c$, let $r'_k$ be the number of barred $k$'s in the rows strictly below row $i+1$. Let $n_k = \max(\mu_k + r_k - \mu_{k-1} - r'_{k-1},0)$, for $k = 2, \dots, c$. What is the significance of $n_k$? Since the row word of $T$ is Yamanouchi, in row $i+1$ of $T$ there are at least $n_k$ barred $k-1$'s, and in row $i$ of $T$ there are at least $n_k$ barred $k$'s. Thus, in row $i+1$ of $T^1$ there must be at least $n_k$ barred $k$'s, and in row $i$ of $T^1$ there must be at least $n_k$ barred $k-1$'s.

For each $2 \leq k \leq c$, define subtableaux $P_k$ and $Q_k$ of $T^1$, where $P_k$ is the subtableau in row $i$ of $T^1$ containing $n_k$ barred $k-1$'s, counting right from the leftmost $k$ in row $i$, and $Q_k$ is the subtableau in row $i+1$ of $T^1$ containing $n_k$ barred $k$'s, counting left from the rightmost $k$ in row $i+1$. Thus, the rightmost box of $P_k$ contains a barred $k-1$ and the leftmost box of $Q_k$ contains a barred $k$. 

We will now form $\tilde{T}$ from $T^1$. We slightly abuse notation to let us communicate the process without requiring messy subscripts.  For each $ k = 2, 3, \dots, c$, do the following independently: Let $n = n_k$, $P = P_k$, and $Q = Q_k$. The subtableau $P$ consists of $n$ blocks of unbarred $k-1$'s, so that each block is to the left of a barred $k-1$. Let  $v_{i}$, $i = 1, \dots, n$, be the number of entries in each block respectively, reading the blocks left to right. Similarly, the subtableau $Q$ consist of $n$ blocks of unbarred $k$'s, so that each block is to the right of a barred $k$. Let $w_i$, $i = 1, \dots, n$, be the number of entries in each block respectively, reading the blocks left to right. Replace the entries in $P$ with $n$ blocks, each containing $v_{i}$ unbarred $k$'s, $i =1,\dots, n$, with each block to the right of a barred $k$. Call the subtableau of entries replacing $P$ in this way $\tilde{Q}_k$. Also, replace the entries in $Q$ with $n$ blocks, each containing $w_{i}$ unbarred $k-1$'s, $i = 1, \dots, n$, with each block to the left of a barred $k-1$. Call the subtableau of entries replacing $Q$ in this way $\tilde{P}_k$. The process works by exchanging $n_k$ barred $k$'s in row $i+1$ with the same amount of barred $k-1$'s in row $i$, while keeping the weights of the affected unbarred entries unchanged (we will check this later). The tableau formed after applying this process independently to each pair of $P_k$, and $Q_k$, $k = 2,\dots, c$ is the tableau $\tilde{T}$, and this completes the process of \emph{reordering barred entries} and we have finished describing the involution on the set of barred tableaux.

\subsection{Applying $\psi_1$, then $\psi_2$, is an involution}
\ble
Let $T$ be a tableau bad in row $i$, column $j$. Denote by $\psi: T \to \tilde{T}$ the process of applying $\psi_1$ then $\psi_2$ to $T$ according to the diagram:
\[
T \overset{\psi_1}{\longrightarrow} T^1 \overset{\psi_2}{\longrightarrow} \tilde{T}.
\]
Then we claim $\psi$ is an involution, that is, $\psi:\tilde{T} \to T$.
\ele

\bpf

As a result of the tail swap (process $\psi_1$), we have the following properties of $\tilde{T}$:
\begin{itemize}
	\item[1.] $\tilde{T}$ is bad in row $i$, column $j$.
	\item[2.] The rows of $\tilde{T}$ weakly decrease, left to right.
	\item[3.] $\tilde{T}$ is of shape $\tilde{\ka} = \ol{R}_{i,i+1} \ka$.
\end{itemize}
Recall that if $T$ obeys condition (C2a), then we had to apply the process of fixing the column ordering. Recall the integers $q$ and $r$ (Definitions \ref{def:q} and \ref{def:r}). There were two cases, Case 1 was for $q \geq r$ and Case 2 was for $q < r$. If $T$ is a tableau that has the property $q \geq r$ then $\tilde{T}$ has the property that $q< r$. Similarly, if $T$ is a tableau that has the property $q < r$ then $\tilde{T}$ has the property that $q\geq r$. Moreover, the process of fixing the column ordering is an involution.

For all $T$, the tail swap when applied to $\tilde{T}$ of shape $\tilde{\ka}$ restores the shape $\ka$. For each $ 2 \leq k \leq c$, the tail swap also sends the subtableaux $\tilde{P_k}$ to row $i$ and the subtableaux $\tilde{Q_k}$ to row $i+1$. Then, the process of reordering barred entries restores the subtableaux ${P_k}$ and ${Q_k}$ to their original locations in $T$.

%
%

\epf

We now give a sequence of lemmas to check that the weights of $T$ and $\tilde{T}$ are almost equal.

\subsection{Weights of $T$ and $\tilde{T}$ are almost equal}
For the purposes of this subsection, let $\rho(\al)$ denote the labelling $\rho^r(\al)$ on $T$, and $\tilde{\rho}(\al)$ denote the labelling $\tilde{\rho}^r(\al)$ on $\tilde{T}$.
\ble
Suppose $T$ obeys condition (C2a) and $q \neq r-1$. Let $\beta = (i,q')$ be a box containing an unbarred entry of $T$ between boxes $(i,j)$ and $(i,q)$, and $\ga = (i+1,r')$ be a box containing an unbarred entry of $T$ between boxes $(i+1,j+1)$ and $(i+1,r)$. Let $\beta'$ and $\gamma'$ be the boxes $(i+1,q')$ and $(i,r')$ respectively; that is, $\beta'$ and $\gamma'$ are the boxes below $\beta$ and above $\ga$ respectively. Then $\eval(T(\be)) = \eval(\tilde{T}(\be'))$ and $\eval(T(\ga)) = \eval(\tilde{T}(\ga'))$.
\ele

\bpf
When $q \neq r-1$, the process of fixing column ordering pairs the entry $T(\beta)$ to the entry $\tilde{T}(\beta')$, and the entry $T(\ga)$ to the entry $\tilde{T}(\ga')$. We claim that the weight of $T(\beta)$ is the same as the weight of $\tilde{T}(\beta')$. Suppose there are $t$ barred $b$'s between box $(i,j)$ and $(i,q')$ of $T$. Then there are $t-1$ barred $c$'s between box $(i+1,j+1)$ and $(i,r')$ of $\tilde{T}$. Since $b = c+1$, we have $a_{b - \rho(\be)_{b}} = a_{c - \rho(\be')_{c}}$ and $a_{b-c(\be)} = a_{c-c(\be')}$. Thus, the weight of   $T(\beta)$ is the same as the weight of $\tilde{T}(\beta')$.  A similar argument shows that the weight of $T(\ga)$ is the same as the weight of $\tilde{T}(\ga')$.
\epf

Note that if $T$ obeys condition (C2a) but $q = r-1$, then the entries in boxes $(i,j)$ to $(i,q)$ and $(i+1,j+1)$ to $(i+1,r)$ of $T$ and $\tilde{T}$ are the same.

\ble
For $2\leq k \leq c$, the weight of $P_k$ is equal to the weight of $\tilde{Q}_k$, and  the weight of $Q_k$ is equal to the weight of $\tilde{P}_k$.
\ele

\bpf
Note that for each $2 \leq k\leq c$ the barred subtableau $P_k$ are paired to the subtableau of entries in $T$ which contain the first $n_k$ barred $k-1$'s in row $i+1$ of T, with respect to the row order. Similarly, the subtableau $Q_k$ is paired to the subtableau of $T$ which contains the last $n_k$ barred $k$'s in row $i$ of $T$.  After barred entries are reordered, the subtableau $\tilde{P}_k$ (which occupies the boxes of $Q_k$) contains the first $n_k$ barred $k-1$'s in row $i+1$ of $\tilde{T}$ and $\tilde{Q}_k$ (which occupies the boxes of $P_k$) contains the last $n_k$ barred $k$'s in row $i$ of $\tilde{T}$. Then, the proof follows in exactly the same manner as the previous proof.
\epf

\ble
The weight of all barred tableaux $T$ of shape $\ka$, bad in row $i$, column $j$, is equal to the weight of all barred tableaux $T$ of shape $\tilde{\ka}$, bad in row $i$, column $j$.
\ele

\bpf
The process of applying $\psi_1$ then $\psi_2$ is almost a weight preserving involution on barred tableaux bad in row $i$, column $j$. This is because the unbarred entries of $T$ which are affected by the process of fixing column ordering (if $T$ obeys condition (C2a)), and entries  in $P_k$ or $Q_k$ for all $2\leq k \leq c$ are paired to entries with corresponding weight in $\tilde{T}$ by the previous two lemmas.  However, the other entries might not have the same weight because of the tail swap. We were unable to find an involution on these entries that would preserve the weight, so we adopt the approach of cancelling paired monomials ocurring in the weight of $T$ and $\tilde{T}$.

We define the \emph{unaffected entries} of $\tilde{T}$ to be the entries of $\tilde{T}$ which are unaffected by the process of fixing column ordering (if $T$ obeys condition (C2a)), and entries strictly to the right of boxes $(i,j-1)$ and $(i+1,j)$ that are not in $\tilde{P}_k$ or $\tilde{Q}_k$ for all $2\leq k \leq c$. Furthermore, we define the \emph{unaffected entries} of $T$ which are the entries of $T$ which are paired with the unaffected entries of $\tilde{T}$ by the tail swap.

For a $1 \leq k \leq c-1$, let $\delta$ denote the subdiagram of $\ka$ which contains the unaffected entries equal to $k$ in row $i$ of $T$. Similarly, let $\epsilon$  denote the subdiagram of $\ka$ which contains the unaffected entries equal to $k$ in row $i+1$ of $T$. In fact, note that these unaffected entries occur after $n_{k+1}$ barred $k$'s in row $i+1$, and before $n_{k}$ barred $k$'s in row $i$ of $T$, with respect to the row order. Let $\tilde{\delta}$ and $\tilde{\epsilon}$ denote the respective subdiagrams of $\tilde{\ka}$ that are paired to $\delta$ and $\epsilon$ by the tail swap. 

Let $M$ be the subtableaux containing barred and unbarred $k$'s of $T$ in the subdiagrams $\delta$ and $\epsilon$. Then the weight of $M$ is
\[
\prod_{\substack{\alpha \in \delta \cup \epsilon \\ T(\alpha) = k \text{ unbarred}}} (a_{T(\alpha)-\rho^r(\alpha)_{T(\alpha)}} - a_{T(\alpha) - c(\alpha)}).
\]
We wish to split this weight into monomials, so recall our definition of monomial tableaux from Section 2. A \emph{monomial subtableau} $N$ is \emph{derived} from $M$ by doing the following: for each unbarred entry in $M$, either add a prime as a superscript of that unbarred entry or do nothing. Then, the weight of $M$ can be expanded as:
\[\begin{split}
\eval(M) &= \sum_N \eval(N)\\
&= \sum_N \prod_{\substack{\alpha \in \delta \cup \epsilon \\ T(\alpha)= k \text{ unbarred,}\\ \text{unprimed}}} (a_{T(\alpha)-\rho^r(\alpha)_{T(\alpha)}})  \prod_{\substack{\alpha \in \delta \cup \epsilon \\ T(\alpha)= k \text{ unbarred,}\\ \text{primed}}}( - a_{|T(\alpha)| - c(\alpha)})
\end{split}\]
summed over all monomial subtableaux $N$ derived from $M$.

Let $N$ be a monomial tableaux derived from $M$. We will find a monomial subtableau $\tilde{N}$ of $\tilde{T}$ in the subdiagrams $\tilde{\delta}\cup\tilde{\epsilon}$ such that the weight of $N$ and $\tilde{N}$ are equal. Let $N'$ denote the sequence of unprimed $k$'s in $M$, listed left to right, first in the subdiagram $\delta$, then in the subdiagram $\epsilon$. Let $\chi$ be the subset of boxes of the subdiagram $\delta \cup \epsilon$ containing the primed entries of $M$. Let $\tilde{\chi}$ be the subset of the subdiagram $\tilde{\delta} \cup \tilde{\epsilon}$ paired to $\chi$ by the process of tail swapping. To form $\tilde{N}$ first fill in the boxes of $\tilde{\chi}$ with primed $k$'s. Then, fill in the boxes of $\tilde{\delta} \cup \tilde{\epsilon}$ not in $\tilde{\chi}$ with unprimed $k$'s such that the sequence of replaced entries, read left to right, first from row $i+1$ and then from row $i$, is equal to $N'$. This forms the monomial tableaux $\tilde{N}$.  Note that these replaced entries occur after $n_{k+1}$ barred $k$'s in row $i+1$ of $\tilde{T}$, and before $n_k$ barred $k$'s in row $i$ of $\tilde{T}$, thus the weight of $N$ and $\tilde{N}$ are equal.
\epf
This completes the proof of Theorem \ref{thm:lr} since we have cancelled out all unwanted summands from (\ref{eqn:weaklr}).

We may express the Littlewood--Richardson polymonials using the following alternative form, which is equivalent to \cite[Theorem 2.1]{Molev:lrp}.
\bco\label{cor:lr}
Let $\la$, $\mu$, $\nu$ be partitions. If $\nu \not\subseteq \mu$, then $\clmn = 0$. If $\mu \subseteq \nu$, we have that
\beq\label{eqn:lr2}
\clmn = \sum_{T} \prod_{\substack{\alpha \in \la\\ T(\alpha) \text{ unbarred}}}(a_{T(\alpha)-\rho^c(\alpha)_{T(\alpha)}} - a_{T(\alpha) - c(\alpha)}),
\eeq
where the sum is taken over reverse $\la$-tableaux $T$ obeying the following. First, the column word $S^c$ of $T$ is Yamanouchi and $S^c:\mu \to \nu$. Secondly, the entries in $T$ strictly decrease down each column.
\eco

\bpf
The difference between Corollary \ref{cor:lr} and Theorem \ref{thm:lr} is the usage of the labelling $\rho^c(\alpha)$ instead of $\rho^r(\alpha)$. The corollary follows from the fact that the entries in $T$ strictly decrease down each column.
\epf

\bre
We may apply the $\nu$-boundedness condition from \cite[Theorem 2.1]{Molev:lrp} to make our formula Graham positive \cite{Graham}.
\ere

\bre\label{rmk:alternativerule}
In Tamvakis \cite{Tamvakis} the Littlewood--Richardson rule for $\clmn$ is given in terms of skew tableaux of shape $\nu/\mu$. It may be possible to use the same idea to obtain a rule to calculate the polynomials $\clmn(a)$ which depends on skew tableaux of shape $\nu/\mu$.
\ere

\subsection{Example}
\bex
We give an example of a bad tableau $T$ which falls in Subcase 1a of the proof.  Let $\mu = (2^2)$, $\kappa= (9^2)$, and $\nu = (4,3,2,1)$. Then the following is a bad tableau appearing in $\kkmn(a)$:
\[
T=\yvc\ybd\young(\fourb\four\four\four32\twob22,3\three\threeb3\threeb31\oneb\oneb)
\]
Let $\alpha = (i,j) = (1,1)$.   Since $L(\alpha) = 4$, we have that $L(\alpha)$ takes $\mu$ to $(2,2,0,1)$, which is not a partition. Thus, the word $L(\al)$ is not Yamanouchi, (P1) holds, and $T$ is bad. We have the entries $b= \ol{4}$, and $c = 3$, in boxes $(1,1)$ and $(2,1)$ respectively. In row 1, box $(1,q) = (1,4)$ is the rightmost box containing a 4, barred or unbarred, and there are $s = 1$ barred $4$'s between box (1,1) and box (1,4). Counting right from box (2,2), we see there is 1 barred $3$ up until box $(2,r) = (2,3)$. There are $x_1 = 3$ unbarred 4's between box $\alpha$ and box $(1,4)$ and $y_1 = 1$ unbarred 3's between box $(2,2)$ and $(2,3)$. We have $r \leq q$, so $T$ obeys condition (C2a), Case 1. We obtain the following tableau $T^{\frac{1}{2}}$ after we fix the column ordering:
\[
T^{\frac{1}{2}} = \yvc\ybd\young(\boldfourb\boldfour\boldthree\boldthreeb32\twob22,3\boldthree\boldthree3\threeb31\oneb\oneb)
\]
where the bold entries are the ones which have been affected. Now,  we see that we have $y_1 = 1$ unbarred 4's, and $x_1 = 3$ unbarred 3's. The entries weakly decrease down column 3. Thus, let $X$ be the subtableau of entries from box $(1,3)$ to $(1,9)$, and $Y$ the subtableau of entries from box $(2,4)$ to $(2,9)$; $X$ and $Y$ are the bold entries in row 1 and 2 respectively of the tableau:
\[
T^{\frac{1}{2}} = \yvc\ybd\young(\fourb4\boldthree\boldthreeb\boldthree\boldtwo\boldtwob\boldtwo\boldtwo,333\boldthree\boldthreeb\boldthree\boldone\boldoneb\boldoneb)
\]
We perform a tailswap on $X$ and $Y$, producing the tableau $T^1$:
\[
T^1=\yvc\ybd\young(\fourb4\boldthree\boldthreeb\boldthree\boldone\boldoneb\boldone,333\boldthree\boldthreeb\boldthree\boldtwo\boldtwob\boldtwo\boldtwo)
\]
$T^1$ has shape $ \ol{R}_{1,2}(9,9) = (8,10) $ and each row of $T^1$ weakly decreases from left to right. We claim that the weight of $\eval(T^1(2,3))$ is equal to the weight of $\eval(T(1,3))$. This is because $3 - \mu_3 = 4 - \rho^r(1,3)$, and $c(2,3) = c(1,3) -1$ We now apply the second process: reordering the barred entries. We calculate $n_2$ and $n_3$. Since there are no entries below row $2$, $r'_1 = r'_2 = 0$. The number of barred $2$'s and $3$'s in rows 1 and 2 are $r_2 = 1$ and $r_3 = 2$ respectively. Then, $n_2 =\max(\mu_2 + r_2 - \mu_1 - r'_1 , 0) = 1$ and $n_3 = \max(\mu_3 + r_3 - \mu_2 - r'_2, 0) = 0$. Thus, we find subtableau $P_2$ and $Q_2$ of $T$ containing one barred 1 and one barred 2 respectively. The subtableau $P_2$ and $Q_2$ consists of the bold entries in row 1 and row 2 respectively of the tableau:
\[ T^1=\yvc\ybd\young(\fourb43\threeb3\boldone\boldoneb\oneb,3333\threeb32\boldtwob\boldtwo\boldtwo),
\]
It is a fact that the row word of $T^1$ is not Yamanouchi but we will fix shortly. We have that $v_1 = 1$ and $w_1 = 2$; these are the numbers of unbarred $1$'s and $2$'s in $P_2$ and $Q_2$ respectively. We obtain $\tilde{P_2}$ and $\tilde{Q_2}$, which are the bold entries in row $2$ and row $1$ respectively of the tableau:
\[ \tilde{T} =\ybd\young(\fourb43\threeb3\boldtwob\boldtwo\oneb,3333\threeb32\boldone\boldone\boldoneb),
\]
which completes the process of reordering barred entries and thus we have formed $\tilde{T}$ from $T$. We claim that the weight of $\tilde{T}(1,7)$ is equal to the weight of $T(2,7)$. The entry $T(2,7)$ has no barred 1's before it. The entry $\tilde{T}(1,7)$ has one barred $2$ before it. The weights are equal since $1 - \rho^r(2,7)_1 = 2 - \tilde{\rho}^r(1,7)_2$ and $c(2,7) = c(1,7) -1$.  

The entries of $\tilde{T}$ not in $\tilde{P}_2$ and $\tilde{Q}_2$ which are also unaffected by fixing the column ordering are marked in bold:
\[
\ybd\young(\fourb4\boldthree\boldthreeb\boldthree\twob2\boldoneb,3333\threeb\boldthree\boldtwo11\oneb)
\]
We will check that the weights of the unaffected entries can be cancelled out later, but first we check that applying the involution to $\tilde{T}$ will restore $T$ (throughout bold entries denote affected entries). We do the tail swap, first obtaining
\[
\tilde{T}^{\frac{1}{2}} = \yvc \ybd\young(\fourb4\boldthree\boldthreeb\boldthree\boldtwo\boldone\boldone\boldoneb,333\boldthree\boldthreeb\boldthree\boldtwob\boldtwo\boldoneb)
\]
from which we fix the column ordering and obtain:
\[
\tilde{T}^{\frac{1}{2}} =\yvc \ybd\young(\fourb\boldfour\boldfour\boldfour3211\oneb,3\boldthree\boldthreeb3\threeb3\twob2\oneb)
\]
Then, we reorder the barred entries, obtaining:
\[
T =  \yvc\ybd\young(\fourb44432\boldtwob\boldtwo\boldtwo,33\threeb3\threeb3\boldone\boldoneb\oneb)
\]
Thus, the process restores $T$ from $\tilde{T}$.

Consider the subtableaux of unaffected entries in $\tilde{T}$, marked in bold, which are the entries not affected by fixing the column ordering and not in $\tilde{P}_2$ or $\tilde{Q}_2$ of $\tilde{T}$:
\beq\label{eqn:m1}
\yvc\ybd\young(\fourb4\boldthree\boldthreeb\boldthree\twob2\boldoneb,3333\threeb\boldthree\boldtwo11\oneb)
\eeq
This tableau originated from the bad tableau $T$ of shape (9,9):
\beq\label{eqn:m2}
\yvc\ybd\young(\fourb444\boldthree\boldtwo\twob22,33\threeb\boldthree\boldthreeb\boldthree1\oneb\boldoneb)
\eeq
with the bold entries of tableau (\ref{eqn:m2}) paired to the bold entries of tableau (\ref{eqn:m1}) by the tail swap. By the previous arguments the weight of the nonbold entries in tableau (\ref{eqn:m2}) and (\ref{eqn:m1}) are equal.

From the bold entries in tableau (\ref{eqn:m2}) we form the following monomial subtableaux, we are only concerned with the bold entries of tableau (\ref{eqn:m2}) equal to $3$ so we omit the rest:
\[
N = \yvc\ybd\young(\hfil\hfil\hfil\hfil3\hfil\hfil\hfil\hfil,\hfil\hfil\hfil\threep\threeb3\hfil\hfil\hfil)
\]
From this tableau, we have the integer sequence $N' = \ol{3}33$ formed by listing the $3$'s left to right and omitting any primed 3's, first from row 1, then row 2. We form the following subtableaux:
\[
\tilde{N} = \yvc\ybd\young(\hfil\hfil\threep33\hfil\hfil\hfil,\hfil\hfil\hfil\hfil\hfil\threeb\hfil\hfil\hfil\hfil)
\]
which again has the sequence $N' = \ol{3}33$ when the unprimed entries of $\tilde{N}$ are listed, omitting primed entries, left to right, starting from row $2$, then row $1$. Since the boxes containing the primed $3$ in $N$ and $\tilde{N}$ have the same content, the weight of the primed entries are equal. We claim that the weight of the unprimed entries are equal. This follows when we consider that $N$ is the monomial subtableau marked in bold inside the following tableau
\[
\yvc\ybd\young(\fourb444\boldthree2\twob22,33\threeb\boldthreep\boldthreeb\boldthree1\oneb\oneb)
\]
and $\tilde{N}$ is the monomial subtableau marked in bold inside the tableau
\[
\yvc\ybd\young(\fourb4\boldthreep\boldthree\boldthree\twob21,3333\threeb\boldthreeb211\oneb)
\]
where the nonbold entries of the above two tableau have equal weight, by previous arguments. Note that the bold entries occur after one barred $3$ in both tableaux, with respect to the row order. Thus, the weights of the unprimed bold entries in both tableaux are equal, and we conclude that the weights of $N$ and $\tilde{N}$ are equal.
\eex

\end{document}